\documentclass[12pt]{elsarticle}

\usepackage{amsmath,amssymb,amsfonts,mathtools,bm}
\usepackage{hyperref}
\usepackage{enumitem}
\usepackage{booktabs}
\usepackage{algorithm}
\usepackage{algorithmic}
\usepackage{microtype}
\usepackage{mathrsfs}
\usepackage{graphicx}
\usepackage{multicol}
\usepackage{multirow}
\usepackage{caption, subcaption}
\usepackage[dvipsnames]{xcolor}

\usepackage[letterpaper,left=2.5cm,right=2.5cm,top=2.5cm,bottom=2.5cm]{geometry}

\allowdisplaybreaks
\hypersetup{colorlinks=true,linkcolor=blue,citecolor=blue,urlcolor=blue,
 pdftitle={High-order fully discrete multi-entropy-stable and bound-preserving schemes for relativistic Euler equations},
 pdfauthor={Linfeng Xu; Kailiang Wu}}

\newtheorem{theorem}{Theorem}[section]
\newtheorem{remark}[theorem]{Remark}
\newtheorem{example}[theorem]{Example}
\newcounter{testcase}[section]
\renewcommand{\thetestcase}{\thesection.\arabic{testcase}}
\newcommand{\testcase}[2]{%
  \refstepcounter{testcase}\par\medskip
  \noindent\textbf{Test \thetestcase: #1}\label{#2}\par\smallskip}

\newcommand{\R}{\mathbb R}
\newcommand{\E}{\mathcal E}
\newcommand{\Qflux}{\mathcal Q}
\newcommand{\Qhat}{\widehat{\mathcal Q}}
\newcommand{\Fhat}{\widehat{\mathbf F}}
\newcommand{\Gset}{\mathcal G}
\newcommand{\Geps}{\mathcal G_\varepsilon}

\newcommand{\dd}{\,d}

\begin{document}

\begin{frontmatter}

\title{Fully discrete multi-entropy-stable bound-preserving high-order schemes for relativistic Euler equations }

\author[a]{Linfeng Xu}
\author[a,b]{Kailiang Wu\corref{cor1}}
\cortext[cor1]{Corresponding author.}

\affiliation[a]{organization={Department of Mathematics, Southern University of Science and Technology}, city={Shenzhen}, postcode={518055}, country={China}}
\affiliation[b]{organization={Shenzhen International Center for Mathematics, Southern University of Science and Technology}, city={Shenzhen}, postcode={518055}, country={China}}

\begin{abstract}
A discrete entropy inequality is the principal nonlinear stability estimate available for systems of conservation laws, and evaluating it presupposes a physically admissible state. So far, however, the two have been secured separately. Entropy-stable schemes are almost always semi-discrete, are built around one selected entropy pair, and take for granted the positivity of density and pressure that makes the entropy well defined in the first place, whereas bound-preserving limiters keep the solution admissible but deliver no entropy estimate. For the special relativistic Euler equations, the two cannot be separated at all, since the conservative-to-primitive map is implicit, and an inadmissible state therefore has no entropy to correct. Here we construct high-order discontinuous Galerkin and finite volume schemes that, to our knowledge, for the first time, are entropy
stable in the fully discrete sense for an arbitrary prescribed finite family of convex entropy pairs, a property we call multi-entropy stability, and are provably admissible wherever an entropy is evaluated. All of this is achieved by a single cellwise projection, and neither conservation nor the design order is lost. The construction rests on relativistic causality, which bounds every characteristic speed by the speed of light. Consequently, the numerical viscosity can be fixed once for all states and all equations of state, and one
two-point building block then serves the whole entropy family. Since only the convexity of the admissible set and this speed bound are used, the same route remains open for related systems. Finally, in computations with four equations of state, the schemes retain high-order accuracy close to vacuum, produce no inadmissible state in strong shocks, near-vacuum shock--vortex interaction or jets with Lorentz factor above $70$, and confirm the monotone decay of every enforced discrete entropy.
\end{abstract}


\end{frontmatter}

\section{Introduction}\label{sec:intro}

Relativistic hydrodynamics (RHD) governs flows whose bulk velocity approaches the speed of light, or whose internal energy becomes comparable with the rest-mass energy. Such conditions are routine in high-energy astrophysics, in extragalactic and magnetized relativistic jets \cite{van1996knots,hughes2002three}, in gamma-ray bursts, pulsar winds, and core-collapse supernovae, and they are approached in the laboratory in inertial-confinement fusion and in laser--plasma interaction. Since the early computations of May and White
\cite{may1966hydrodynamic} and Wilson \cite{wilson1972numerical}, numerical simulation has been the principal tool in these regimes, and a substantial body of shock-capturing methodology has been built for the special
relativistic Euler equations, surveyed in
\cite{marti2003numerical,font2008numerical,rezzolla2013relativistic,marti2015grid}. In the strongly relativistic and near-vacuum regimes of interest, however, a high-order scheme repeatedly produces states outside the admissible set, and robustness rather than resolution has remained the limiting factor.

That admissible set, however, is not described explicitly in the variables the scheme evolves. Conservative and primitive variables are coupled through the Lorentz factor and the specific enthalpy, and for a general equation of state (EOS), the inverse map is available only through a scalar root-finding problem
\cite{WuTang2017ApJS,ryu2006equation,cai2024provably}. If a numerical state has a nonpositive rest-mass density or pressure, or a velocity at or beyond the speed of light, the safeguarded iteration (see, for example, \ref{app:primitive}) has no root to converge to. The eigenvalues of the flux Jacobian and the specific entropy are then both undefined at that quadrature node, so the time step cannot be completed at all. Consequently, simulations of ultra-relativistic jets or near-vacuum regions break down within a few steps, in our experience, during the initial transient, unless admissibility is enforced by construction.

Admissibility alone, however, does not single out the right solution. The relativistic Euler system is nonlinear, its solutions develop discontinuities from smooth data, and weak solutions are not unique, so an entropy inequality is needed to select the physically relevant one \cite{dafermos2000hyperbolic}. The same is true after discretization, where a discrete entropy inequality is one of the few nonlinear stability estimates
available for a system of conservation laws \cite{tadmor2003entropy}, and it excludes convergence to a spurious weak solution.

These two requirements interact, because every entropy evaluation passes through the primitive recovery. A limiter that computes an entropy correction and afterward repairs positivity is therefore not well defined, since an inadmissible state has no entropy to correct. The obstruction is more severe here than for the Newtonian Euler equations, where the pressure is at least an explicit function of the conservative variables and the specific entropy is available in closed form. In the relativistic case, admissibility itself has to be characterized indirectly in the conservative variables, and every entropy evaluation requires a scalar solve that is well posed only on the admissible set. The entropy correction is consequently defined only once admissibility has been established, and positivity must accordingly be enforced first.

We shall therefore require a high-order discretization that
\begin{itemize}
	\item a discrete entropy inequality holds for the \emph{fully discrete} update, and not for the semi-discrete method of lines alone,
	\item the same conservative update control several convex entropy pairs, and
	\item every entropy evaluation is carried out on a state whose admissibility has already been established.
\end{itemize}
These requirements are met below for discontinuous Galerkin (DG) and finite volume discretizations of the
special relativistic Euler equations with a general causal EOS of the form \eqref{eq:gEOS}, by one cellwise convex projection.

\subsection{Bound-preserving schemes and entropy-stable schemes}

Physical-constraint-preserving, bound-preserving, and invariant-domain-preserving methods have developed mainly through scaling limiters and flux correction. The scaling limiter of Zhang and Shu \cite{zhang2010,zhang2010b} first establishes admissibility of the updated cell average and then contracts the high-order polynomial toward that average. It has been carried to the compressible Euler equations, to special and general relativistic hydrodynamics, and to relativistic magnetohydrodynamics \cite{QinShu2016,WuTangM3AS,WuTang2017ApJS,WuShu2018,WuShu2019,wu2017design},
and a survey of this line is given in \cite{zhang2011b}. Flux-corrected transport and the more recent convex limiting instead blend a high-order update with a robust low-order one
\cite{boris1973flux,book1975flux,boris1976flux,zalesak1979fully,zalesak2012design,
Hu2013,WuTang2015,guermond2017invariant,guermond2018second,guermond2019invariant,
abgrall2024bound,abgrall2025bound,abgrall2025novel,xu2026gql}, and geometric quasilinearization provides a general algebraic route to such bounds \cite{wu2023geometric}, with recent reviews given in \cite{kuzmin2024property,wu2025high}. These techniques are by now standard in high-order relativistic computation, but they provide no entropy estimate.

Entropy-conservative and entropy-stable discretizations form a parallel line of development, from the work of Harten, Osher, Tadmor, and Jiang and Shu to high-order finite volume, DG, summation-by-parts, and flux-differencing methods \cite{harten1976finite,crandall1980monotone,osher1984riemann,osher1988convergence,
tadmor1987numerical,tadmor2003entropy,jiang1994cell,lefloch2002fully, fjordholm2012arbitrarily,fjordholm2013eno,hiltebrand2014entropy,chen2017entropy}.
Formulations are available for many fluid systems, including RHD and relativistic magnetohydrodynamics
\cite{ismail2009affordable,chandrashekar2013kinetic,liu2018entropy,duan2019high, bhoriya2020entropy,duan2021entropy,wu2020entropy,duan2020high,xu2024high}. In this literature, however, the entropy estimate is usually semi-discrete, it is stated for a single entropy pair, and the entropy variables and entropy fluxes are evaluated under the assumption that the density and the pressure are positive. For a relativistic model with a general EOS that last assumption is not a convenience, since the entropy is computable only when it holds.

Fully discrete entropy control has been pursued through space--time DG \cite{hiltebrand2014entropy}, entropy correction terms \cite{abgrall2022reinterpretation}, relaxation Runge--Kutta methods \cite{ranocha2020relaxation,xu2024high}, algebraic flux correction \cite{kuzmin2022limiter}, and, most recently, a limiter-based construction for hyperbolic conservation laws \cite{liu2026limiter}. Positivity has also been
combined with a bound on the specific entropy inside a single DG limiter for chemically reacting flows \cite{ching2024positivity,ching2024positivity2}. For RHD, a minimum principle on the specific entropy and the invariant-region-preserving schemes built on it were established in \cite{WuMEP2021}, with a local version for high-order schemes in \cite{cui2025local}. An entropy-based indicator has been used in \cite{guercilena2017entropy} to control the local reconstruction order. A pointwise bound on the specific entropy is, however, a different statement from a discrete entropy inequality, and none of these constructions handles a family of entropy pairs, or the requirement that admissibility be established before any entropy is evaluated.

To the best of our knowledge, no scheme for the relativistic Euler equations is at present simultaneously fully discrete, entropy stable for more than one entropy pair, and provably admissible at every point at which an entropy is evaluated. The schemes constructed below have all three properties.

\subsection{Multi-entropy stability and the proposed construction}

Our starting points are the weak-to-strong scaling mechanism and the entropy--positivity--oscillation (EPO) framework of Wu \cite{wu2026epo}, which isolates the convex geometry shared by entropy limiting, bound preservation, and oscillation control. We do not reprove that geometry: the radial projection, the SSP convex-combination argument, and the extension to Cartesian meshes are taken from \cite{wu2026epo}. 
A similar entropy limiting approach was also recently presented in the independent work of Liu et al.~\cite{liu2026limiter}. 
The present paper develops the relativistic content that the EPO framework leaves open: the invariant
domain and its convexity in conservative variables, the uniform light-cone bound that fixes the numerical viscosity for all admissible data, the transfer of entropy admissibility from the thermodynamic pair to an entire prescribed family, the ordering forced by the implicit primitive recovery, and the fully discrete schemes and computations that follow from them.

Fix an arbitrary finite family $\{(\E^{(\ell)},\Qflux^{(\ell)})\}_{\ell=1}^M$ of known convex entropy pairs. We
call a fully discrete update \emph{multi-entropy stable} when a single conservative state satisfies the local and global discrete entropy inequalities for every $\ell=1,\ldots,M$. Each pair defines a feasible interval on a common scaling ray, and the intersection of the $M$ intervals is used once. Applying $M$ unrelated single-entropy corrections in sequence would not suffice, because each correction may undo the effect of its predecessors.

All directional characteristic speeds of the relativistic Euler system lie in $[-1,1]$ for any causal EOS, so the Lax--Friedrichs viscosity may be set to $\alpha=1$ uniformly, with no wave-speed estimate and with no dependence on the state or on the EOS. Under the standard hypothesis that the directional Riemann problem has a self-similar entropy solution lying in the admissible set, the resulting two-point average is the average of that solution over the light cone, hence admissible, and it inherits one entropy inequality per pair. A Gauss--Lobatto endpoint decomposition then writes one forward Euler step of the high-order scheme as a convex combination of nodal states and two such two-point averages, which turns the two-point inequalities into a weak entropy budget, one per pair, for the updated cell average. The candidate polynomial is next contracted along one ray anchored at that average. The density and the concave admissibility functional are enforced first, and the entropies are evaluated only on the shortened, admissible part of the ray. Their common feasible interval, intersected with an optional convex oscillation-suppressing interval, then gives the final scaling factor, and since the anchor is the cell average, the scheme remains conservative. Finally, SSP multistep methods are convex combinations of forward Euler steps \cite{Gottlieb2009}, so both the budgets and the projection carry over to the fully discrete schemes.

\subsection{Contributions and organization}

The principal contributions are the following.
\begin{enumerate}[label=(\roman*),leftmargin=2em]
\item We construct high-order DG and finite volume schemes for the special relativistic Euler equations with a general causal EOS. The schemes are fully discrete, keep every quadrature state in the relativistic invariant domain, and satisfy local and global discrete entropy inequalities for an arbitrary prescribed finite family of convex entropy pairs. The SSP multistep realization uses a single projection per time step.
\item We show that every directional characteristic speed of the system lies in $[-1,1]$ for any causal EOS. The Lax--Friedrichs coefficient can therefore be fixed at $\alpha=1$, independently of the state and of the EOS, and the resulting two-point average supplies one entropy inequality per pair. Moreover, only the thermodynamic entropy pair has to be assumed admissible across a discontinuity, because the generating function satisfies $\mathcal H'(S)>0$, admissibility for every other member of the family follows. No entropy differentiation or inversion of the entropy-variable map is needed.
\item Explicit weak entropy budgets for the relativistic Euler system are derived in one and two dimensions for the same high-order candidate average. Adding an entropy pair only adds scalar entropy evaluations inside one synchronized root-finding procedure, whereas the conservative update and the number of primitive recoveries per trial radius remain unchanged.
\item The weak-to-strong step is realized by a projection that enforces admissibility before entropy. Consequently, no primitive recovery and no entropy evaluation is ever attempted outside the admissible set. The projection
preserves the cell average, restores nodal admissibility, enforces all $M$ entropy budgets, and accommodates a convex oscillation-control module on the same ray.
\item Smooth tests confirm the design order to within $10^{-5}$ of vacuum. An ultra-relativistic Riemann problem isolates the effect of using several entropy pairs and shows that they remove an overshoot left by a single pair behind the contact discontinuity. Further one- and two-dimensional tests with four EOSs assess strong shocks, shock heating, near vacuum, and Lorentz factors above $70$. In every reported entropy history, the discrete total entropy, corrected for the boundary flux on non-periodic domains, decreases monotonically.
\end{enumerate}

Section~\ref{sec:preliminaries} states the model, the causal EOS class, the invariant domain, and the entropy family. Section~\ref{sec:twopoint} develops the light-cone bound and the two-point analysis specific to RHD.
Section~\ref{sec:weakbudgets} derives the weak multi-entropy budgets for high-order candidate averages. Section~\ref{sec:epo} constructs the projection, and Section~\ref{sec:time} embeds it in fully discrete SSP time integrators. Section~\ref{sec:numerics} reports the numerical assessment, and Section~\ref{sec:conclusion} concludes. Primitive recovery and the entropy formulas for the four EOSs are presented in the appendices.

\section{The relativistic Euler equations, admissible states, and entropy pairs}\label{sec:preliminaries}

\subsection{Conservative variables and a causal equation of state}

We consider the $d$-dimensional special relativistic Euler equations
\begin{equation}\label{eq:pde}
\partial_t \mathbf U + \sum_{i=1}^d \partial_{x_i}\mathbf F_i(\mathbf U)=0,
\end{equation}
with
\begin{equation*}
\mathbf U=\begin{pmatrix}D\\ \mathbf m\\ E\end{pmatrix}
=\begin{pmatrix}\rho\gamma\\ \rho h\gamma^2\mathbf v\\ \rho h\gamma^2-p\end{pmatrix},
\qquad
\mathbf F_i(\mathbf U)=\begin{pmatrix}Dv_i\\ v_i\mathbf m+p\mathbf e_i\\ m_i\end{pmatrix}.
\end{equation*}
Here $D,\ {\mathbf m}$, and $E$ are the mass density, the momentum density vector, and the energy density. These are called the conservative variables. The associated primitive variables are the rest-mass density $\rho$, the pressure $p$ and the velocity $\mathbf v=(v_1,\dots,v_d)^\top$ with $|\mathbf v|<1$. The Lorentz factor is $\gamma=(1-|\mathbf v|^2)^{-1/2}$, $e$ denotes the specific internal energy and $h=1+e+p/\rho$ the specific enthalpy. Throughout, the speed of light is normalized to one.

The relativistic Euler system \eqref{eq:pde} is closed by an equation of state (EOS), which relates the specific enthalpy $h$ to pressure $p$ and rest-mass density $\rho$. 
In this work, we  assume that $h$ only depends on $\theta:=\frac{p}{\rho}$ and takes the form 
\begin{equation} \label{eq:gEOS}
h = h(\theta) = 1  + e (\theta) + \theta.
\end{equation}
To ensure relativistic causality (local sound speed $0<c_s(\theta)<1$), the EOS must satisfy the inequality \cite{WuTang2017ApJS}:
\begin{equation*}
    h'(\theta)>\frac{h(\theta)}{h(\theta)-\theta},\qquad \theta>0.
\end{equation*}

\subsection{Admissible states in conservative variables}

The physically admissible states are naturally expressed in primitive variables as
\begin{equation*}
\Gset_{\mathrm{phys}}
:=
\{\mathbf U=(D,\mathbf m,E)^\top: \rho(\mathbf U)>0,\ p(\mathbf U)>0,\ |\mathbf v(\mathbf U)|<1\}.
\end{equation*}
However, because $\rho$, $p$, and $\mathbf{v}$ are not available as explicit functions of the conservative variables, the following equivalent characterization directly in conservative variables is required \cite{WuTang2015}:
\begin{equation*}
\Gset:=\{\mathbf U=(D,\mathbf m,E)^\top:D>0,\ \mathfrak q(\mathbf U):=E-\sqrt{D^2+|\mathbf m|^2}>0\}.
\end{equation*}
The map $(D,\mathbf m)\mapsto\sqrt{D^2+|\mathbf m|^2}$ is convex, so
$\mathfrak q$ is concave and $\Gset$ is open and convex. Consequently, convex
averages of admissible states are admissible, and each physical constraint can
be enforced along a ray issuing from an admissible mean.

\subsection{A prescribed finite family of convex entropy pairs}

The relativistic Euler system admits an infinite family of convex entropy pairs, constructed in \cite{cui2025local} for EOS \eqref{eq:gEOS}. We recall only what is needed below. Define the specific entropy
\begin{equation*}
S(\mathbf U)=-\ln\rho+\int_1^\theta\frac{e'(\xi)}{\xi}\,\dd\xi.
\end{equation*}
For any smooth $\mathcal H:\R\to\R$ satisfying
\begin{equation}\label{eq:Hcond}
\mathcal H'(S)>0,\qquad
\mathcal H'(S)-\bigl(1+e'(\theta)\bigr)\mathcal H''(S)>0,
\end{equation}
the entropy pair is
\begin{equation*}
\E(\mathbf U)=-D\mathcal H(S(\mathbf U)),\qquad
\Qflux_i(\mathbf U)=-Dv_i\mathcal H(S(\mathbf U)).
\end{equation*}
The canonical and physically most important choice is $\mathcal H(S)=S$, which yields the thermodynamic entropy pair
\begin{equation}\label{eq:canonicalpair}
\E(\mathbf U)=-DS(\mathbf U),
\qquad
\Qflux_i(\mathbf U)=-Dv_iS(\mathbf U).
\end{equation}

We fix an arbitrary but finite family
\begin{equation}\label{eq:entropy-family-M}
\mathfrak E_M:=\{(\E^{(\ell)},\Qflux_1^{(\ell)},\ldots,\Qflux_d^{(\ell)})\}_{\ell=1}^M,
\end{equation}
where each $\E^{(\ell)}\in C^2(\Gset)\cap C(\overline{\Gset})$ is strictly
convex on $\Gset$. The analysis is uniform in $M$, since an additional pair changes
only scalar entropy evaluations and leaves the conservative candidate update
unchanged.

\subsection{High-order nodal representation and local discrete entropy}

In one dimension, let $I_j=[x_{j-1/2},x_{j+1/2}]$ and $x_j^{(\nu)}$, $\nu=1,\ldots,L$, be the nodes of a quadrature rule with positive weights that is exact for the local polynomial degree. Given a candidate DG polynomial $\mathbf U_h|_{I_j}$, we write
$
\mathbf U_{j,\nu}:=\mathbf U_h(x_j^{(\nu)}),\;
\mathbf U_j := (\mathbf U_{j,1},\dots,\mathbf U_{j,L}),
$
and its cell average
$
\bar{\mathbf U}_j := \frac1{|I_j|}\int_{I_j} \mathbf U_h(x)\,\dd x.
$
With quadrature weights $\omega_\nu>0$ satisfying $\sum_{\nu=1}^L \omega_\nu=1$, exactness gives
\begin{equation*}
\bar{\mathbf U}_j = \sum_{\nu=1}^L \omega_\nu \mathbf U_{j,\nu}.
\end{equation*}
For each entropy pair, define
\begin{equation}\label{eq:localdiscE}
\E_j^{(\ell)}(\mathbf U_j):=\sum_{\nu=1}^L\omega_\nu\E^{(\ell)}(\mathbf U_{j,\nu}).
\end{equation}
Jensen's inequality gives
$\E^{(\ell)}(\bar{\mathbf U}_j)\le\E_j^{(\ell)}(\mathbf U_j)$. Here, the left-hand side is a cell-average quantity, and the right-hand side is the nodal quantity that the method controls. Our aim is to design a conservative fully discrete scheme that keeps all quadrature states admissible while satisfying a local entropy inequality for every member of \eqref{eq:entropy-family-M}.

\section{The first-order building block: a light-speed bound and a two-point entropy inequality}\label{sec:twopoint}

The EPO construction of \cite{wu2026epo} requires a two-point average that is admissible and that satisfies a two-point entropy inequality. As we now show, relativistic causality supplies such an average with a single viscosity coefficient valid for all admissible data, and that same average then serves the whole prescribed entropy family.

\subsection{A uniform bound on the directional characteristic speeds}\label{subsec:alpha1}

Let $\mathbf n\in\R^d$ be a unit vector. Set the directional flux, the normal velocity, and the squared tangential speed as
\begin{equation*}
\mathbf F_{\mathbf n}(\mathbf U):=\sum_{i=1}^d n_i\mathbf F_i(\mathbf U),\qquad
v_n:=\mathbf v\cdot\mathbf n,\qquad v_t^2:=|\mathbf v|^2-v_n^2.
\end{equation*}
The projected system has the repeated speed $\lambda_0=v_n$ and acoustic speeds
\begin{equation*}
\lambda_\pm=\frac{v_n(1-c_s^2)\pm \dfrac{c_s}{\gamma}
\sqrt{1-v_n^2-c_s^2v_t^2}}{1-c_s^2|\mathbf v|^2};
\end{equation*}
see \cite{WuTang2017ApJS}.
Set $s=c_s^2\in(0,1)$, $w=v_n$, $\tau=v_t^2$, and write
\begin{equation*}
\lambda_\pm=\frac{w(1-s)\pm A}{B},\qquad
B:=1-s(w^2+\tau)>0,
\end{equation*}
\begin{equation*}
A:=\sqrt{s(1-w^2-\tau)(1-w^2-s\tau)}.
\end{equation*}
A direct factorization gives
\begin{align*}
(B-w(1-s))^2-A^2&=(1-s)(1-s|\mathbf v|^2)(1-w)^2\ge0,\\
(B+w(1-s))^2-A^2&=(1-s)(1-s|\mathbf v|^2)(1+w)^2\ge0.
\end{align*}
Moreover,
\begin{align*}
B-w(1-s)&\ge(1-w)(1-s)\ge0,\\
B+w(1-s)&\ge(1+w)(1-s)\ge0.
\end{align*}
Hence $|w(1-s)\pm A|\le B$, and therefore
\begin{equation*}
\max\{|v_n|,|\lambda_-|,|\lambda_+|\}\le1
\quad\text{for every }\mathbf U\in\Gset\text{ and every unit }\mathbf n.
\end{equation*}
The Lax--Friedrichs viscosity can consequently be fixed at $\alpha=1$.
All subsequent two-point statements are then free of any wave-speed estimate, and the CFL restriction depends only on the mesh and on the quadrature rule.

\subsection{From canonical entropy admissibility to the entire family}

For a unit vector $\mathbf n$ and two states $\mathbf U_L,\mathbf U_R\in\Gset$, we define the directional local Lax--Friedrichs average with $\alpha=1$ by
\begin{equation*}
H_{\mathbf n}(\mathbf U_L,\mathbf U_R)
:=\frac{\mathbf U_L+\mathbf U_R}{2}
-\frac{\mathbf F_{\mathbf n}(\mathbf U_R)-\mathbf F_{\mathbf n}(\mathbf U_L)}{2}.
\end{equation*}
As is standard, we assume that the Riemann problem in the direction $\mathbf n$ has a self-similar solution that remains in $\Gset$ and satisfies the canonical entropy inequality. Only the canonical pair requires an entropy assumption. The reason is that, across a jump moving with speed $\sigma$, mass conservation makes $m:=D(v_n-\sigma)$ continuous, and every pair generated by $\mathcal H$ obeys
\begin{equation*}
[\Qflux_{\mathbf n}]-\sigma[\E]=-m[\mathcal H(S)].
\end{equation*}
For $\mathcal H(S)=S$, entropy admissibility amounts to $m[S]\ge0$. Since $\mathcal H'(S)>0$, the jump $[\mathcal H(S)]$ has the same sign as $[S]$, and the same Riemann solution is therefore admissible for every member of \eqref{eq:entropy-family-M}.

The light-cone bound confines the Riemann fan to $|\xi|\le1$. Integrating the conservation law over this cone yields
\begin{equation*}
H_{\mathbf n}(\mathbf U_L,\mathbf U_R)
=\frac12\int_{-1}^1\mathcal R_{\mathbf n}(\xi;\mathbf U_L,\mathbf U_R)\,\dd\xi.
\end{equation*}
Since $\Gset$ is convex, this average is again admissible,
$H_{\mathbf n}(\mathbf U_L,\mathbf U_R)\in\Gset.$
For each entropy pair $\ell\in\{1,\ldots,M\}$, the Riemann entropy inequality and Jensen's inequality imply
\begin{equation}\label{eq:twopointE}
\begin{aligned}
\E^{(\ell)}\!\left(H_{\mathbf n}(\mathbf U_L,\mathbf U_R)\right)
\le{}\frac{\E^{(\ell)}(\mathbf U_L)+\E^{(\ell)}(\mathbf U_R)}{2}
-\frac{\Qflux_{\mathbf n}^{(\ell)}(\mathbf U_R)
-\Qflux_{\mathbf n}^{(\ell)}(\mathbf U_L)}{2}.
\end{aligned}
\end{equation}

Consequently, a single state $H_{\mathbf n}$ and one viscosity coefficient support all $M$ inequalities, without differentiating the entropy or inverting an entropy-variable map.

\section{Weak multi-entropy budgets for high-order cell averages}\label{sec:weakbudgets}

In this section, the two-point building block of Section~\ref{sec:twopoint} is extended to the cell average produced by one forward-Euler step of a high-order scheme. This leads to a pair of statements that we refer to as \emph{weak budgets}: the updated mean $\bar{\mathbf U}_j^\star$ lies in $\Gset$, and it satisfies one weak entropy inequality for each entropy pair.

\subsection{One-dimensional candidate states}

Let $\omega_\nu$ be the Gauss--Lobatto weights satisfying
$\omega_1=\omega_L$ and $\sum_{\nu=1}^L\omega_\nu=1.$
At time $t^n$, the LF flux with the universal coefficient $\alpha=1$ is
\begin{equation}\label{eq:LFflux1D}
\Fhat_{j+1/2}^n
:=\frac{\mathbf F(\mathbf U_{j,L}^n)+\mathbf F(\mathbf U_{j+1,1}^n)}{2}
-\frac{\mathbf U_{j+1,1}^n-\mathbf U_{j,L}^n}{2}.
\end{equation}
The forward-Euler candidate average is
\begin{equation} \label{eq:candidateavg}
\bar{\mathbf U}_j^\star=\bar{\mathbf U}_j^n
-\lambda(\Fhat_{j+1/2}^n-\Fhat_{j-1/2}^n),\qquad
\lambda:=\frac{\Delta t}{\Delta x}.
\end{equation}

\subsubsection{Weak bound preservation}
Substituting the quadrature identity and \eqref{eq:LFflux1D} into \eqref{eq:candidateavg} leads to
\begin{equation} \label{eq:geom-decomp}
\bar{\mathbf U}_j^\star
=\sum_{\nu=2}^{L-1}\omega_\nu\mathbf U_{j,\nu}^n
+(\omega_1-\lambda)\mathbf U_{j,1}^n
+(\omega_L-\lambda)\mathbf U_{j,L}^n
+\lambda H(\mathbf U_{j-1,L}^n,\mathbf U_{j+1,1}^n)
+\lambda H(\mathbf U_{j,1}^n,\mathbf U_{j,L}^n).
\end{equation}
If the input quadrature states $\mathbf U_{j,\nu}^n$ are admissible and
\begin{equation}\label{eq:CFL}
    0<\lambda\le\omega_1,
\end{equation}
then \eqref{eq:geom-decomp} is a convex combination of admissible states, so $\bar{\mathbf U}_j^\star\in\Gset$.

\subsubsection{One weak budget per entropy pair}
Define the numerical entropy flux as
\begin{equation}\label{eq:Qhat1D}
\Qhat_{j+1/2}^{n,(\ell)}
:=\frac{\Qflux^{(\ell)}(\mathbf U_{j,L}^n)+\Qflux^{(\ell)}(\mathbf U_{j+1,1}^n)}{2}
-\frac{\E^{(\ell)}(\mathbf U_{j+1,1}^n)-\E^{(\ell)}(\mathbf U_{j,L}^n)}{2}.
\end{equation}
Then for each $\ell\in\{1,\ldots,M\}$, applying \eqref{eq:twopointE} to the two $H$ states in
\eqref{eq:geom-decomp} yields
\begin{equation}\label{eq:weakfamily}
\E^{(\ell)}(\bar{\mathbf U}_j^\star)
\le B_j^{n,(\ell)}
:=\sum_{\nu=1}^L\omega_\nu\E^{(\ell)}(\mathbf U_{j,\nu}^n)
-\lambda\bigl(\Qhat_{j+1/2}^{n,(\ell)}-\Qhat_{j-1/2}^{n,(\ell)}\bigr).
\end{equation}
On a periodic mesh,
\begin{equation*}
\sum_jB_j^{n,(\ell)}=\sum_j\E_j^{(\ell)}(\mathbf U_j^n),
\qquad \ell=1,\ldots,M.
\end{equation*}
The decomposition \eqref{eq:geom-decomp} is independent of $\ell$: all entropy pairs share the same conservative update and the same admissible mean, and only the budgets differ.

\begin{remark}[Dissipation and the two two-point averages]
Fixing $\alpha=1$ avoids wave-speed estimation but may be more dissipative than a local Rusanov speed in slowly moving regions. A sharper flux can replace \eqref{eq:LFflux1D} only if it supplies a compatible two-point inequality. 
The two $H$ states in \eqref{eq:geom-decomp} are dictated by the flux decomposition itself: one couples the values of the outer neighbours, the other the endpoints of the current cell. Their inequalities recombine into the telescoping interface fluxes \eqref{eq:Qhat1D}.
\end{remark}

\subsection{Two-dimensional Cartesian meshes}

For a Cartesian cell $K_{ij}=I_i\times J_j$, we adopt an $L$-point Gauss--Lobatto rule in the normal direction and a $Q$-point Gauss rule in the transverse direction. The local DG polynomial has the line-separable average
\begin{equation*}
\bar{\mathbf U}_{ij}^n
=\ \kappa_1\sum_{q=1}^Q\widetilde\omega_q
\sum_{\mu=1}^L\widehat\omega_\mu
p_{ij}^n(\widehat x_i^{(\mu)},\widetilde y_j^{(q)})
+\kappa_2\sum_{q=1}^Q\widetilde\omega_q
\sum_{\mu=1}^L\widehat\omega_\mu
p_{ij}^n(\widetilde x_i^{(q)},\widehat y_j^{(\mu)}).
\end{equation*}
The standard weights are
\begin{equation*}
\kappa_1=\frac{\alpha_1/\Delta x}{\alpha_1/\Delta x+\alpha_2/\Delta y},\qquad
\kappa_2=\frac{\alpha_2/\Delta y}{\alpha_1/\Delta x+\alpha_2/\Delta y},
\end{equation*}
where $\alpha_1, \alpha_2$ are the maximum absolute characteristic speeds in the $x$- and $y$-directions, and we simply take $\alpha_1=\alpha_2=1$ by Section~\ref{subsec:alpha1}.

Define the one-dimensional slice averages
$$\bar{\mathbf U}_{ij,q}^{x,n} := \sum_{\mu=1}^L \widehat{\omega}_\mu p_{ij}^n(\widehat{x}_i^{(\mu)}, \widetilde{y}_j^{(q)}),
\qquad \bar{\mathbf U}_{ij,q}^{y,n} := \sum_{\mu=1}^L \widehat{\omega}_\mu p_{ij}^n(\widetilde{x}_i^{(q)}, \widehat{y}_j^{(\mu)}),$$
so that
$$\bar{\mathbf U}_{ij}^n = \kappa_1 \sum_{q=1}^Q \widetilde{\omega}_q \bar{\mathbf U}_{ij,q}^{x,n} + \kappa_2 \sum_{q=1}^Q \widetilde{\omega}_q \bar{\mathbf U}_{ij,q}^{y,n}.$$

\subsubsection{Weak bound preservation}
Let $\lambda_x := {\Delta t}/{\Delta x}, \lambda_y := {\Delta t}/{\Delta y}$. The forward-Euler candidate average is 
\begin{equation*}
\bar{\mathbf U}_{ij}^\star = \bar{\mathbf U}_{ij}^n
- \lambda_x \sum_{q=1}^Q \widetilde{\omega}_q \left( \Fhat_{1,i+1/2,q}^n - \Fhat_{1,i-1/2,q}^n \right)
- \lambda_y \sum_{q=1}^Q \widetilde{\omega}_q \left( \Fhat_{2,q,j+1/2}^n - \Fhat_{2,q,j-1/2}^n \right),
\end{equation*}
where the directional LF fluxes are evaluated at the edge Gauss nodes:
$$\Fhat_{1,i+1/2,q}^n := \Fhat_1\left(p_{ij}^n(x_{i+1/2}, \widetilde{y}_j^{(q)}), p_{i+1,j}^n(x_{i+1/2}, \widetilde{y}_j^{(q)})\right),$$
$$\Fhat_{2,q,j+1/2}^n := \Fhat_2\left(p_{ij}^n(\widetilde{x}_i^{(q)}, y_{j+1/2}), p_{i,j+1}^n(\widetilde{x}_i^{(q)}, y_{j+1/2})\right).$$
Introducing the directional candidate slice averages
\begin{align*}
{\mathbf A}_{ij,q}^{x,\star} := \bar{\mathbf U}_{ij,q}^{x,n} - \frac{\lambda_x}{\kappa_1} \left( \Fhat_{1,i+1/2,q}^n - \Fhat_{1,i-1/2,q}^n \right),\quad
{\mathbf A}_{ij,q}^{y,\star} := \bar{\mathbf U}_{ij,q}^{y,n} - \frac{\lambda_y}{\kappa_2} \left( \Fhat_{2,q,j+1/2}^n - \Fhat_{2,q,j-1/2}^n \right),
\end{align*}
the forward-Euler update admits the exact convex representation
\begin{equation}\label{eq:2Dconvexrep}
\bar{\mathbf U}_{ij}^\star = \kappa_1 \sum_{q=1}^Q \widetilde{\omega}_q {\mathbf A}_{ij,q}^{x,\star} + \kappa_2 \sum_{q=1}^Q \widetilde{\omega}_q {\mathbf A}_{ij,q}^{y,\star}.
\end{equation}
Each $\mathbf A$ is a one-dimensional update. Thus, admissibility follows when
$\frac{\alpha_1\Delta t}{\kappa_1\Delta x}\le\widehat\omega_1,
\;
\frac{\alpha_2\Delta t}{\kappa_2\Delta y}\le\widehat\omega_1,
$
or, with the standard $\kappa$ values,
\begin{equation*}
\left(\frac{\alpha_1}{\Delta x}+\frac{\alpha_2}{\Delta y}\right)\Delta t
\le\widehat\omega_1.
\end{equation*}

\subsubsection{One weak budget per entropy pair}
For each entropy pair, define the directional edge entropy fluxes
\begin{equation*}
\Qhat_{1,i+1/2,q}^{(\ell),n} := \frac{\Qflux_1^{(\ell)}({\mathbf U}_{i+1/2,q}^{-,n}) + \Qflux_1^{(\ell)}({\mathbf U}_{i+1/2,q}^{+,n})}{2}
- \frac{\alpha_1}{2} \left( \E^{(\ell)}({\mathbf U}_{i+1/2,q}^{+,n}) - \E^{(\ell)}({\mathbf U}_{i+1/2,q}^{-,n}) \right),
\end{equation*}
\begin{equation*}
\Qhat_{2,q,j+1/2}^{(\ell),n} := \frac{\Qflux_2^{(\ell)}({\mathbf U}_{q,j+1/2}^{-,n}) + \Qflux_2^{(\ell)}({\mathbf U}_{q,j+1/2}^{+,n})}{2}
- \frac{\alpha_2}{2} \left( \E^{(\ell)}({\mathbf U}_{q,j+1/2}^{+,n}) - \E^{(\ell)}({\mathbf U}_{q,j+1/2}^{-,n}) \right),
\end{equation*}
where ${\mathbf U}_{i+1/2,q}^{\pm,n}$ and ${\mathbf U}_{q,j+1/2}^{\pm,n}$ are the one-sided states at the edge Gauss nodes. The local discrete entropy on $K_{ij}$ is defined as
\begin{equation} \label{eq:2DdiscE}
\E_{ij}^{(\ell)}({\mathbf U}_{ij}) := \kappa_1 \sum_{q=1}^Q \widetilde{\omega}_q \sum_{\mu=1}^L \widehat{\omega}_\mu \E^{(\ell)}\bigl({\mathbf U}_{ij}(\widehat{x}_i^{(\mu)}, \widetilde{y}_j^{(q)})\bigr)
+ \kappa_2 \sum_{q=1}^Q \widetilde{\omega}_q \sum_{\mu=1}^L \widehat{\omega}_\mu \E^{(\ell)}\bigl({\mathbf U}_{ij}(\widetilde{x}_i^{(q)}, \widehat{y}_j^{(\mu)})\bigr).
\end{equation}
The slice decomposition \eqref{eq:2Dconvexrep} reduces the two-dimensional entropy estimate to one-dimensional applications of \eqref{eq:weakfamily}. By convexity of the entropy function $\E^{(\ell)}$, we obtain
\begin{align}
\E^{(\ell)}(\bar{\mathbf U}_{ij}^\star) \leq B_{ij}^{(\ell),n} := \E_{ij}^{(\ell)}({\mathbf U}_{ij}^n) &- \lambda_x \sum_{q=1}^Q \widetilde{\omega}_q \left( \Qhat_{1,i+1/2,q}^{(\ell),n} - \Qhat_{1,i-1/2,q}^{(\ell),n} \right)\notag\\
&- \lambda_y \sum_{q=1}^Q \widetilde{\omega}_q \left( \Qhat_{2,q,j+1/2}^{(\ell),n} - \Qhat_{2,q,j-1/2}^{(\ell),n} \right).
\label{eq:2DweakE}
\end{align}

\section{The multi-entropy projection}\label{sec:epo}

The weak budgets of Section~\ref{sec:weakbudgets} constrain only the cell average. We now convert them into the nodal statements required by the high-order scheme, using a single conservative scaling ray anchored at that mean. Physical admissibility is settled first, since neither the primitive recovery nor the entropy evaluation is defined outside $\Gset$. Each prescribed entropy pair then contributes one feasible interval on the shortened, physically valid part of the ray. The intersection of these intervals is combined with an optional oscillation-control interval. The underlying radial geometry is established in \cite{wu2026epo}.

\subsection{Candidate arrays and a common scaling ray}

The weak budgets derived above refer to a forward Euler candidate average. The local EPO analysis uses only the data
\begin{equation}\label{eq:abstractweak}
\bar{\mathbf U}_j^\star\in \Gset,
\qquad
\E^{(\ell)}(\bar{\mathbf U}_j^\star)\le B_j^{(\ell)},
\qquad \ell=1,\dots,M,
\end{equation}
together with a candidate nodal array
$
\mathbf U_j^\star = (\mathbf U_{j,1}^\star,\dots,\mathbf U_{j,L}^\star)
$
whose weighted average is $\bar{\mathbf U}_j^\star$. 
For every $\theta\in[0,1]$, define the scaling ray anchored at the candidate average by
\begin{equation*}
\mathcal S_{\bar{\mathbf U}_j^\star}(\theta;\mathbf U_j^\star)
:=
\Bigl(
\bar{\mathbf U}_j^\star+\theta(\mathbf U_{j,1}^\star-\bar{\mathbf U}_j^\star),\dots,
\bar{\mathbf U}_j^\star+\theta(\mathbf U_{j,L}^\star-\bar{\mathbf U}_j^\star)
\Bigr).
\end{equation*}
At $\theta=0$ all nodal states collapse to the cell average, and at $\theta=1$ we recover the original candidate nodal array. All points along the ray have the same weighted average $\bar{\mathbf U}_j^\star$, so conservation holds for every $\theta\in[0,1]$.

\subsection{The positivity radius}

The positivity module computes the largest common scaling factor that keeps every nodal state admissible. Since $\Gset$ is convex, each nodal ray can leave $\Gset$ at most once.

For each node, define
\begin{equation*}
\mathbf U_{j,\nu}(\theta)
:=
\bar{\mathbf U}_j^\star+\theta(\mathbf U_{j,\nu}^\star-\bar{\mathbf U}_j^\star),
\qquad 0\le \theta\le1.
\end{equation*}
We split admissibility into the two conditions $D>0$ and $\mathfrak q>0$.

\begin{itemize}[leftmargin=*]
\item Density radius: 
Set a cell-local threshold $0<\varepsilon_D\leq\bar D_j^\star$. In practice we take $\varepsilon_D=\min\{\varepsilon,\bar D_j^\star\}$ with a global $\varepsilon>0$, which is a valid choice for every $\bar{\mathbf U}_j^\star\in\Gset$.
The density radius is
\begin{equation}\label{eq:thetaD}
\theta_j^D
:=
\min\left\{
1,\
\min_{\nu:\,D_{j,\nu}^\star<\bar D_j^\star}
\frac{\bar D_j^\star-\varepsilon_D}{\bar D_j^\star-D_{j,\nu}^\star}
\right\}.
\end{equation}

\item $\mathfrak q$-radius:
Since $\mathfrak q$ is concave on $\Gset$, for each node $\nu$, Jensen's inequality gives 
\begin{equation} \label{q-radius-Jensen}
\mathfrak q(\mathbf U_{j,\nu}(\theta))
\geq
(1-\theta) \mathfrak q(\bar{\mathbf U}_j^\star) + \theta \mathfrak q(\mathbf U_{j,\nu}^\star).
\end{equation}
Take a cell-local threshold $0<\varepsilon_q\leq\mathfrak q(\bar{\mathbf U}_j^\star)$, for example $\varepsilon_q=\min\{\varepsilon,\mathfrak q(\bar{\mathbf U}_j^\star)\}$.
If $\mathfrak q(\mathbf U_{j,\nu}^\star)\ge \mathfrak q(\bar{\mathbf U}_j^\star)$, then the right-hand side of \eqref{q-radius-Jensen} is at least $\mathfrak q(\bar{\mathbf U}_j^\star)\geq\varepsilon_q$. Otherwise, we require the linear lower bound in \eqref{q-radius-Jensen} to be no smaller than $\varepsilon_q$, so that $\mathfrak q(\mathbf U_{j,\nu}(\theta))\geq\varepsilon_q$. 
The global $\mathfrak q$-radius is then
\begin{equation}\label{eq:thetaq}
\theta_j^{\mathfrak q}
:=
\min\left\{
1,\
\min_{\nu:\,\mathfrak q(\mathbf U_{j,\nu}^\star)<\mathfrak q(\bar{\mathbf U}_j^\star)}
\frac{\mathfrak q(\bar{\mathbf U}_j^\star)-\varepsilon_\mathfrak q}{\mathfrak q(\bar{\mathbf U}_j^\star)-\mathfrak q(\mathbf U_{j,\nu}^\star)}
\right\}.
\end{equation}
\end{itemize}
The density is affine, and for $\mathfrak q$, concavity guarantees the bound \eqref{eq:thetaq} without solving nonlinear equations. Hence,
\begin{equation*}
\theta_j^P := \min\{\theta_j^D,\theta_j^\mathfrak q\}
\end{equation*}
is well defined. Every smaller radius remains feasible, so for each node $\nu$ and $0\le \theta\le \theta_j^P$,
\begin{equation*}
\mathbf U_{j,\nu}(\theta)\in \Geps
:=\{\mathbf U=(D,\mathbf m,E)^\top:\ D\ge \varepsilon_D,\ \mathfrak q(\mathbf U)\ge \varepsilon_q\}\subset\Gset.
\end{equation*}

The positivity-limited array is
\begin{equation*}
\mathbf U_j^P:=\mathcal S_{\bar{\mathbf U}_j^\star}(\theta_j^P;\mathbf U_j^\star).
\end{equation*}

\subsection{One entropy radius per prescribed pair}

Fix one entropy pair $(\E^{(\ell)},\Qflux^{(\ell)})$ in the family. Because $\mathbf U_j^P$ lies entirely in $\Geps\subset \Gset$, the entropy is well defined along the shortened ray
$
\mathcal S_{\bar{\mathbf U}_j^\star}(\vartheta;\mathbf U_j^P),
\;
0\le \vartheta\le 1.
$
A relativistic entropy cannot be evaluated without a primitive-variable recovery (\ref{app:primitive}), and that recovery is well-posed only for states in $\Gset$. The unshortened candidate ray generally leaves $\Gset$, where the entropy is not defined. The positivity module must therefore act first. Once $\mathbf U_j^P$ lies in $\Gset$, the entropy profile on the shortened ray is well defined. 

Define the entropy profile on the shortened ray as
\begin{equation}\label{eq:Psij}
\Psi_j^{(\ell)}(\vartheta)
:=
\sum_{\nu=1}^L \omega_\nu
\E^{(\ell)}\Bigl(
\bar{\mathbf U}_j^\star + \vartheta(\mathbf U_{j,\nu}^P-\bar{\mathbf U}_j^\star)
\Bigr).
\end{equation}
The profile $\Psi_j^{(\ell)}$ is convex as it is a convex combination of the convex functions $\E^{(\ell)}$ evaluated along affine paths. Since all nodal states on the ray have the same weighted mean $\bar{\mathbf U}_j^\star$,
Jensen's inequality yields $\Psi_j^{(\ell)}(\vartheta)\ge\E^{(\ell)}(\bar{\mathbf U}_j^\star) =\Psi_j^{(\ell)}(0)$. A convex function attaining its minimum at the left endpoint is nondecreasing on $[0,1]$. Therefore, we have
\begin{equation*}
\Psi_j^{(\ell)}(\vartheta)\ge \Psi_j^{(\ell)}(0)=\E^{(\ell)}(\bar{\mathbf U}_j^\star)
\quad\text{for all }0\le \vartheta\le1.
\end{equation*}

This monotonicity turns a weak entropy budget into a strong one. The strong entropy $\Psi_j^{(\ell)}(1)=\E_j^{(\ell)}(\mathbf U_j^P)$ may violate the budget $B_j^{(\ell)}$, but the weak budget ensures $\Psi_j^{(\ell)}(0)\le B_j^{(\ell)}$. Because the profile is convex and nondecreasing, the set of $\vartheta$ for which $\Psi_j^{(\ell)}(\vartheta)\le B_j^{(\ell)}$ is a closed interval containing $0$. Compressing the nodal values toward the mean reduces the entropy, so a suitable $\vartheta$ can always be found.

For each $\ell=1,\dots,M$, define the individual feasible endpoint
\begin{equation*}
\vartheta_j^{(\ell)}
:=
\sup\Bigl\{\vartheta\in[0,1]:\ \Psi_j^{(\ell)}(\vartheta)\le B_j^{(\ell)}\Bigr\}.
\end{equation*}
Their intersection is represented by the single multi-entropy radius
\begin{equation*}
\vartheta_j^{\mathrm{ME}}
:=\min_{1\le\ell\le M}\vartheta_j^{(\ell)}
=\sup\Bigl\{\vartheta\in[0,1]:
\Psi_j^{(\ell)}(\vartheta)\le B_j^{(\ell)}\ \text{for every }\ell\Bigr\}.
\end{equation*}
Equivalently, it is the zero-sublevel endpoint of the common feasibility profile
\begin{equation}\label{eq:familyprofile}
\Phi_j(\vartheta):=\max_{1\le\ell\le M}
\bigl(\Psi_j^{(\ell)}(\vartheta)-B_j^{(\ell)}\bigr),
\end{equation}
which is convex and nondecreasing. The family radius, expressed in the original candidate coordinates, is therefore
\begin{equation*}
\theta_j^{PE}:=\theta_j^P\vartheta_j^{\mathrm{ME}}.
\end{equation*}
The weak budget makes $\vartheta=0$ feasible for every pair, while the profile analysis makes each feasible set a nonempty closed interval. Hence $\vartheta_j^{\mathrm{ME}}\in[0,1]$, the radius $\theta_j^{PE}$ is well defined, and the common projected state satisfies
\begin{equation*}
\E_j^{(\ell)}\Bigl(\mathcal S_{\bar{\mathbf U}_j^\star}(\theta_j^{PE};\mathbf U_j^\star)\Bigr)
\le B_j^{(\ell)},
\qquad
\ell=1,\dots,M.
\end{equation*}

\subsection{Optional convex oscillation module}

Oscillations may still be present in the solution produced by a high-order FV/DG scheme after the admissibility and entropy radii have been imposed. In order to use the same radial scaling mechanism for their suppression, we adopt the convex oscillation-suppressing (COS) operator of \cite{cao2026cos}.

For each cell $I_j$ and candidate average $\bar{\mathbf U}_j^\star$, assume we have a nonempty closed convex set $\mathcal O_j^\star$ of nodal arrays that contains the constant array with all entries equal to the candidate mean:
$(\bar{\mathbf U}_j^\star,\dots,\bar{\mathbf U}_j^\star) \in \mathcal O_j^\star.$
Then we define the oscillation radius as
\begin{equation*}
\theta_j^O
:=
\sup\Bigl\{\theta\in[0,1]: \mathcal S_{\bar{\mathbf U}_j^\star}(\theta;\mathbf U_j^\star)\in \mathcal O_j^\star\Bigr\}.
\end{equation*}
Because $\mathcal O_j^\star$ is closed and the ray map is continuous, the set of feasible $\theta$ is closed. Hence, $\theta_j^O$ is well defined and belongs to $[0,1]$.
The convexity of $\mathcal O_j^\star$ then implies
\begin{equation*}
\mathcal S_{\bar{\mathbf U}_j^\star}(\theta;\mathbf U_j^\star)\in \mathcal O_j^\star
\qquad\text{for all }0\le \theta\le \theta_j^O.
\end{equation*}

\subsubsection{An explicit entropy-induced COS module}\label{subsec:cos}

The COS metric requires one strictly convex entropy, but it does not determine which entropy inequalities are enforced. We use the canonical thermodynamic entropy $\E^{\mathrm{can}}=-DS$ only for this metric. The multi-entropy projection continues to enforce all $M$ prescribed budgets. Freezing the canonical entropy Hessian at an admissible cell average yields a natural local quadratic metric for the COS operator of \cite{cao2026cos}. The resulting indicator is invariant under uniform rescaling of the state and requires no problem-dependent tuning beyond the constant $C_k$.

Let $\mathbf U_h^n|_{I_j}\in \mathbb{P}^k(I_j)^{d+2}$ be the DG polynomial in cell $I_j$, and let $\bar{\mathbf U}_j^n\in\Gset$ be its cell average. Define
\begin{equation*}
\mathbf H_j^n:=\nabla^2 \E^{\mathrm{can}}(\bar{\mathbf U}_j^n).
\end{equation*}
Because $\E^{\mathrm{can}}\in C^2(\Gset)$ is strictly convex, $\mathbf H_j^n$ is symmetric positive definite. For any vector-valued polynomial $\mathbf Z$ on an interval $K$, define the local entropy-weighted norm
\begin{equation*}
\|\mathbf Z\|_{L_{\mathrm{can},j}^2(K)}^2
:=
\frac1{|K|}\int_K \mathbf Z(x)^\top \mathbf H_j^n \mathbf Z(x)\,\dd x.
\end{equation*}
If $\mathbf U_j^n$ is viewed as a polynomial and extended to the neighbouring cells $I_{j\pm1}$, we set
\begin{equation*}
\Theta_j^\pm
:=
\|\mathbf U_{j\pm1}^n-\bar{\mathbf U}_j^n\|_{L_{\mathrm{can},j}^2(I_{j\pm1})}^2
\end{equation*}
and
\begin{equation}\label{eq:COSsigma}
\sigma_j^\pm
:=
\begin{cases}
C_k
\dfrac{
\|\mathbf U_{j\pm1}^n-\mathbf U_j^n\|_{L_{\mathrm{can},j}^2(I_{j\pm1})}^2
}{
\Theta_j^\pm
},
& \Theta_j^\pm\ge \varepsilon_j,\\[2ex]
0,
& \Theta_j^\pm<\varepsilon_j,
\end{cases}
\qquad
\widehat\sigma_j:=\sigma_j^-+\sigma_j^+,
\end{equation}
where $C_k>0$ depends only on the polynomial degree and
$
\varepsilon_j:=10^{-12}\max\{|\E^{\mathrm{can}}(\bar{\mathbf U}_j^n)|,1\}.
$
The normalization in \eqref{eq:COSsigma} makes $\widehat\sigma_j$ dimensionless and insensitive to uniform scaling of the state.

To convert the smoothness indicator into a damping factor, we choose a local speed bound
\begin{equation}\label{eq:COSbeta}
\beta_j
:=
\max_{\ell=j-1,j,j+1}\max_{\nu=1,\dots,L}
\rho\!\left(\frac{\partial \mathbf F}{\partial \mathbf U}(\mathbf U_{\ell,\nu}^n)\right).
\end{equation}
From the light-cone analysis in Section~\ref{subsec:alpha1}, every eigenvalue of the Jacobian lies in $[-1,1]$, so $0\le \beta_j\le 1.$
One may therefore use the local value \eqref{eq:COSbeta} or simply set $\beta_j=1$ in every cell. The canonical COS coefficient is
\begin{equation}\label{eq:lambdaCOS}
\lambda_j^{\mathrm{COS}}
:=
\exp\!\left(
-\beta_j\frac{\Delta t}{\Delta x_j}\widehat\sigma_j
\right),
\qquad
0<\lambda_j^{\mathrm{COS}}\le 1.
\end{equation}
The corresponding cellwise COS polynomial is the convex scaling
\begin{equation}\label{eq:COSpoly}
\mathbf U_{j}^{\mathrm{COS}}(x)
=
(1-\lambda_j^{\mathrm{COS}})\bar{\mathbf U}_j^n
+
\lambda_j^{\mathrm{COS}}\,\mathbf U_j^n(x).
\end{equation}

For every cell $I_j$, the operator \eqref{eq:COSpoly} has the following properties.
\begin{enumerate}[label=(\roman*),leftmargin=2em]
\item It preserves the cell average:
$
\frac1{|I_j|}\int_{I_j}\mathbf U_j^{\mathrm{COS}}(x)\,\dd x = \bar{\mathbf U}_j^n.
$
\item The set 
$
\mathcal O_j^{\mathrm{COS},n}
:=
\Bigl\{
\mathcal S_{\bar{\mathbf U}_j^n}(\theta;\mathbf U_j^n):\
0\le \theta\le \lambda_j^{\mathrm{COS}}
\Bigr\}
$
is closed and convex, contains the constant state $\bar{\mathbf U}_j^n$, and yields the oscillation radius
$\theta_j^O=\lambda_j^{\mathrm{COS}}$.
\end{enumerate}

\subsubsection{Entropy inheritance under COS}

The COS operation preserves the entropy budget. For each $\ell\in\{1,\dots,M\}$, let $B_j^{(\ell)}$ be any scalar budget satisfying $\E_j^{(\ell)}(\mathbf U_j^n)\le B_j^{(\ell)}.$ Define the prodile
\[
\widehat{\Phi}_j(\theta):=\E_j^{(\ell)}\bigl(\mathcal S_{\bar{\mathbf U}_j^n}(\theta;\mathbf U_j^n)\bigr),
\qquad 0\le \theta\le 1,
\]
which is similar to \eqref{eq:Psij}. The convexity of $\E^{(\ell)}$ and Jensen's inequality then show that $\widehat{\Phi}_j$ is convex and nondecreasing. Consequently, the COS output \eqref{eq:COSpoly} satisfies the entropy budget
$\E_j^{(\ell)}(\mathbf U_j^{\mathrm{COS}})\le B_j^{(\ell)}.$

\begin{remark}[Modal implementation]
If the DG polynomial is written in a modal basis
\[
\mathbf U_j^n(x)=\sum_{r=0}^k \widehat{\mathbf U}_j^{(r)}\phi_j^{(r)}(x),
\qquad
\phi_j^{(0)}\equiv1,
\]
then \eqref{eq:COSpoly} is equivalent to
\[
\widehat{\mathbf U}_j^{(0),\mathrm{new}}=\widehat{\mathbf U}_j^{(0)},
\qquad
\widehat{\mathbf U}_j^{(r),\mathrm{new}}=\lambda_j^{\mathrm{COS}}\widehat{\mathbf U}_j^{(r)},
\quad r\ge 1.
\]
Thus, the COS step leaves the cell average unchanged and multiplies every higher mode by the same scalar factor. 
\end{remark}

\begin{remark}[Cartesian multidimensional COS formula]
On a Cartesian tensor-product cell $K_{ij}=I_i\times J_j$, the same construction is applied to the full tensor-product polynomial $\mathbf U_{ij}^n(x,y)$ about the cell average $\bar{\mathbf U}_{ij}^n$. One defines
\[
\widehat\sigma_{ij}
=
\sigma_{ij}^{x,-}+\sigma_{ij}^{x,+}+\sigma_{ij}^{y,-}+\sigma_{ij}^{y,+},
\]
where, for example,
\[
\sigma_{ij}^{x,+}
=
\begin{cases}
C_k
\dfrac{\|\mathbf U_{i+1,j}^n-\mathbf U_{ij}^n\|_{L_{\mathrm{can},ij}^2(K_{i+1,j})}^2}
{\|\mathbf U_{i+1,j}^n-\bar{\mathbf U}_{ij}^n\|_{L_{\mathrm{can},ij}^2(K_{i+1,j})}^2},
&
\|\mathbf U_{i+1,j}^n-\bar{\mathbf U}_{ij}^n\|_{L_{\mathrm{can},ij}^2(K_{i+1,j})}^2\ge \varepsilon_{ij},\\[2ex]
0,&\text{otherwise},
\end{cases}
\]
and analogously in the other three directions. Then
\[
\lambda_{ij}^{\mathrm{COS}}
=
\exp\!\left(
-\beta_{ij}\Bigl(\frac{\Delta t}{\Delta x_i}+\frac{\Delta t}{\Delta y_j}\Bigr)\widehat\sigma_{ij}
\right),
\quad
\mathbf U_{ij}^{\mathrm{COS}}(x,y)
=
(1-\lambda_{ij}^{\mathrm{COS}})\bar{\mathbf U}_{ij}^n
+
\lambda_{ij}^{\mathrm{COS}}\mathbf U_{ij}^n(x,y).
\]
\end{remark}

\begin{remark}[Local COS activation by an entropy criterion] \label{remark: local COS}
To reduce numerical dissipation and computational cost, the COS procedure may be applied selectively to interfaces that are potentially near shocks, using the Lax-type entropy criterion in \cite{cao2026cos}. Let $\lambda^+$ and $\lambda^-$ denote the largest and smallest eigenvalues of the flux Jacobian. The interface $x_{j+\frac12}$ is marked as a potential shock if the following condition holds, either in the case $\star = +$ or $\star = -$
\[
\lambda^\star(\bar{\mathbf U}_j^n)-\lambda^\star(\bar{\mathbf U}_{j+1}^n)
>
\delta\Bigl(|\lambda^\star(\bar{\mathbf U}_j^n)|+|\lambda^\star(\bar{\mathbf U}_{j+1}^n)|\Bigr),
\]
with $\delta=0.05$. The local variant then replaces the condition $\Theta_j^\pm\ge\varepsilon_j$ in \eqref{eq:COSsigma} by the indicator $\chi_{j\pm\frac12}=1$, where $\chi=1$ on marked interfaces and $0$ otherwise. All other formulas remain unchanged. The multidimensional case is analogous.
\end{remark}

\subsection{The full limiter and the global strong entropy inequality}

We take the minimum of the entropy and COS radii,
\begin{equation*}
\theta_j^\text{EPO}:=\min\{\theta_j^{PE},\theta_j^O\},
\end{equation*}
and define the final limited nodal array as
\begin{equation*}
\widetilde{\mathbf U}_j
:=
\mathcal S_{\bar{\mathbf U}_j^\star}(\theta_j^\text{EPO};\mathbf U_j^\star).
\end{equation*}

If the candidate average satisfies \eqref{eq:abstractweak} and the oscillation set is convex, the minimum of the positivity, multi-entropy, and oscillation radii guarantees four cellwise properties:
\begin{enumerate}[label=(\roman*),leftmargin=2em]
\item mean preservation:
$
\sum_{\nu=1}^L \omega_\nu \widetilde{\mathbf U}_{j,\nu} = \bar{\mathbf U}_j^\star;
$
\item nodal admissibility:
$
\widetilde{\mathbf U}_{j,\nu}\in \Geps\subset \Gset,
\qquad \nu=1,\dots,L;
$
\item strong entropy inequalities for the entire family:
\begin{equation}\label{eq:localstrongfamily}
\E_j^{(\ell)}(\widetilde{\mathbf U}_j)\le B_j^{(\ell)},
\qquad \ell=1,\dots,M;
\end{equation}
\item membership in the oscillation set:
$
\widetilde{\mathbf U}_j\in \mathcal O_j^\star.
$
\end{enumerate}

Summing the local strong entropy inequalities \eqref{eq:localstrongfamily} over all cells gives a global strong entropy inequality. Assume the weak budgets satisfy the compatibility relation
\begin{equation*}
\sum_j B_j^{(\ell)}\le \mathfrak B^{(\ell)}
\qquad \text{for each }\ell=1,\dots,M.
\end{equation*}
Summing the local strong entropy inequalities~\eqref{eq:localstrongfamily} over all cells and inserting the weak budget~\eqref{eq:weakfamily} gives
\[
\sum_{j=1}^N \Delta x\,\E_j^{(\ell)}(\widetilde{\mathbf U}_j)
\le
\sum_{j=1}^N \Delta x\,\E_j^{(\ell)}(\mathbf U_j^n)
- \Delta t \bigl( \Qhat_{N+1/2}^n - \Qhat_{1/2}^n \bigr),
\]
where the boundary numerical entropy fluxes \(\Qhat_{1/2}^n,\Qhat_{N+1/2}^n\) are computed from~\eqref{eq:Qhat1D} using the prescribed exterior states.
Define the boundary-corrected discrete total entropy
\begin{equation*}
\E_{\mathrm{corr}}^{(\ell),\,n}
:=
\Delta x \sum_{j=1}^N \E_j^{(\ell)}({\mathbf U}_j^{n})
+ \Delta t \sum_{k=0}^{n-1} \bigl( \Qhat_{N+1/2}^{k} - \Qhat_{1/2}^{k} \bigr).
\end{equation*}
Then \(\E_{\mathrm{corr}}^{(\ell),\,n+1}\le \E_{\mathrm{corr}}^{(\ell),\,n}\), which gives a rigorous entropy-stability diagnostic for non-periodic simulations. In the periodic case, the boundary fluxes cancel, and this reduces to 
\[
\sum_j \Delta x\,\E_j^{(\ell)}(\widetilde{\mathbf U}_j)\le \sum_j \Delta x\,\E_j^{(\ell)}(\mathbf U_j^n).
\]

\subsection{Two-dimensional case}
The two-dimensional case follows from the one-dimensional results, because the convex decomposition of the cell average and of the entropy budget reduces the line-separable nodal array to one-dimensional slices.
Let ${\mathbf U}_{ij}$ be the candidate nodal array on the line-separable quadrature set, with mean $\bar{\mathbf U}_{ij}^\star$ and budgets $B_{ij}^{(\ell)}$ from \eqref{eq:2DweakE}. Applying the same radial analysis cell by cell with the discrete entropy \eqref{eq:2DdiscE} yields the limited state ${\mathbf U}_{ij}^{\mathrm{EPO}}$ with the following properties:
\begin{enumerate}[label=(\roman*),leftmargin=2em]
    \item every nodal value of ${\mathbf U}_{ij}^\text{EPO}$ belongs to $\Geps\subset\Gset$;
    \item the local strong entropy inequality
    $\E_{ij}^{(\ell)}({\mathbf U}_{ij}^\text{EPO}) \leq B_{ij}^{(\ell)}$
    holds for every entropy pair;
    \item the limited state lies in the convex oscillation-suppressing set.
\end{enumerate}
If the budgets satisfy a global compatibility condition, summing the local inequalities gives a global strong entropy inequality.

For a non-periodic domain, summing the local strong entropy inequalities over all cells leaves boundary terms, just as in the one-dimensional case. With the notation of the weak entropy budget \eqref{eq:2DweakE},
\[
\begin{aligned}
\sum_{i,j} \mathcal E_{ij}^{(\ell)}(\mathbf U_{ij}^\text{EPO})
\le
\sum_{i,j} \mathcal E_{ij}^{(\ell)}(\mathbf U_{ij}^{n})
&- \lambda_x \sum_{j,q} \widetilde \omega_q
\Bigl( \widehat{\mathcal Q}_{1,N_x+1/2,j,q}^{(\ell),n} - \widehat{\mathcal Q}_{1,1/2,j,q}^{(\ell),n} \Bigr) \\
&- \lambda_y \sum_{i,q} \widetilde \omega_q
\Bigl( \widehat{\mathcal Q}_{2,i,q,N_y+1/2}^{(\ell),n} - \widehat{\mathcal Q}_{2,i,q,1/2}^{(\ell),n} \Bigr),
\end{aligned}
\]
where \(N_x\) and \(N_y\) are the numbers of cells in the \(x\)- and \(y\)-directions, and the boundary entropy fluxes are evaluated at the edge Gauss points. The boundary-corrected discrete total entropy
\begin{equation*}
\begin{aligned}
\E_{\mathrm{corr}}^{(\ell),\,n}
:=
\Delta x \Delta y \sum_{i,j} \mathcal E_{ij}^{(\ell)}(\mathbf U_{ij}^{n})
+
\Delta t \sum_{k=0}^{n-1}
\bigg[
& \Delta y \sum_{j,q} \widetilde \omega_q
\Bigl( \widehat{\mathcal Q}_{1,N_x+1/2,j,q}^{(\ell),k} - \widehat{\mathcal Q}_{1,1/2,j,q}^{(\ell),k} \Bigr)
+\\
& \Delta x \sum_{i,q} \widetilde \omega_q
\Bigl( \widehat{\mathcal Q}_{2,i,q,N_y+1/2}^{(\ell),k} - \widehat{\mathcal Q}_{2,i,q,1/2}^{(\ell),k} \Bigr)
\bigg]
\end{aligned}
\end{equation*}
then satisfies \(\E_{\mathrm{corr}}^{(\ell),\,n+1}\le \E_{\mathrm{corr}}^{(\ell),\,n}\). In the periodic case, this reduces to \(\sum_{i,j}\E_{ij}^{(\ell)}(\mathbf U_{ij}^\text{EPO})\le\sum_{i,j}\E_{ij}^{(\ell)}(\mathbf U_{ij}^{n})\).

\subsection{Algorithm}

Algorithm~\ref{alg:epo} presents the one-dimensional procedure. Every step has a directional analog in two and three dimensions.

\begin{algorithm}[h]
\caption{One-step relativistic EPO projection in cell $I_j$}
\label{alg:epo}
\begin{algorithmic}[1]
\STATE Start from a candidate DG or finite volume state array $\mathbf U_j^\star$ and its candidate average $\bar{\mathbf U}_j^\star$.
\STATE Compute the weak geometric and weak entropy budgets for $\bar{\mathbf U}_j^\star$.
\STATE Compute the density radius $\theta_j^D$ from \eqref{eq:thetaD} and the $\mathfrak q$-radius $\theta_j^\mathfrak q$ from \eqref{q-radius-Jensen}.
\STATE Set $\theta_j^P=\min\{\theta_j^D,\theta_j^\mathfrak q\}$ and form the positivity-limited array $\mathbf U_j^P$.
\STATE Evaluate the common feasibility profile $\Phi_j(\vartheta)$ in \eqref{eq:familyprofile}. At each trial radius, recover the $L$ primitive states
once and evaluate all $M$ entropies.
\STATE If $\Phi_j(1)\le0$, set $\vartheta_j^{\mathrm{ME}}=1$. Otherwise, locate the right endpoint of the interval on which $\Phi_j(\vartheta)\le0$.
\STATE Set $\theta_j^{PE}=\theta_j^P\vartheta_j^{\mathrm{ME}}$.
\STATE If the canonical COS module is used, compute $\widehat\sigma_j$ from \eqref{eq:COSsigma}, $\lambda_j^{\mathrm{COS}}$ from \eqref{eq:lambdaCOS}, and set $\theta_j^O=\lambda_j^{\mathrm{COS}}$. Otherwise, compute $\theta_j^O$ from the chosen oscillation set.
\STATE Set $\theta_j^\text{EPO}=\min\{\theta_j^{PE},\theta_j^O\}$ and output $\widetilde{\mathbf U}_j=\mathcal S_{\bar{\mathbf U}_j^\star}(\theta_j^\text{EPO};\mathbf U_j^\star)$.
\end{algorithmic}
\end{algorithm}

\begin{remark}[Choice of the entropy family]
Enforcing several entropy pairs is optional. With $M=1$ and the canonical pair \eqref{eq:canonicalpair}, the limiter reduces to a single-entropy correction. With $M>1$ the same candidate state is projected once, using the minimum of the radii, so no additional limiting passes and no ordering choices are involved. In the experiments of Section~\ref{sec:numerics} we take $\mathcal H_r(S)=(e^{r(S-S_0)}-1)/r$ with three values of $r$, which is permissible whenever $r$ satisfies the convexity condition~\eqref{eq:Hcond}.
\end{remark}

\section{Fully discrete high-order schemes under SSP time stepping}\label{sec:time}

The forward Euler step already provides the admissible cell average and the $M$ weak budgets required by Section~\ref{sec:epo}. Fully discrete high-order time integration is built by writing an SSP method as a convex combination of such forward Euler blocks. The general SSP convexity arguments are given in \cite[Section~7]{wu2026epo}. Since the stagewise SSP Runge--Kutta (RK) method is only first-order in time once the entropy module is active stage by stage \cite{wu2026epo}, we consider an explicit SSP multistep method
\begin{equation*}
\mathbf U^{n+1,\star}=\sum_{r=0}^{k-1}\alpha_r
\left(\mathbf U^{n-r}+\beta_r\Delta tL(\mathbf U^{n-r})\right),
\quad \alpha_r,\beta_r\ge0,\quad \sum_r\alpha_r=1,
\end{equation*}
with SSP coefficient $\mathcal C_{\mathrm{MS}}:=\min_{\beta_r>0}1/\beta_r$. If
$\Delta t\le\mathcal C_{\mathrm{MS}}\Delta t_{\mathrm{FE}}$, the candidate average is admissible, and for every entropy pair
\begin{equation*}
\E^{(\ell)}(\bar{\mathbf U}_j^{n+1,\star})
\le B_j^{n+1,(\ell)}
:=\sum_{r=0}^{k-1}\alpha_r\Bigl(\E_j^{(\ell)}(\mathbf U_j^{n-r})
-\beta_r\lambda(\Qhat_{j+1/2}^{n-r,(\ell)}-\Qhat_{j-1/2}^{n-r,(\ell)})\Bigr).
\end{equation*}
A single projection to the candidate state ${\mathbf U}^{n+1,\star}$ produces the updated solution $\mathbf U^{n+1}$, which preserves the average, keeps all quadrature states in $\Geps$, and enforces the local strong entropy inequality for every member of the family.

We use the following explicit third-order SSP multistep method in Section~\ref{sec:numerics}.
\begin{example}[The third-order SSP multistep scheme]\label{ex:ssp-ms3}
The third-order SSP multistep method reads
\begin{equation*}
\mathbf U^{n+1,\star}
=\frac{16}{27}\left(\mathbf U^n+3\Delta tL(\mathbf U^n)\right)
+\frac{11}{27}\left(\mathbf U^{n-3}+\frac{12}{11}\Delta tL(\mathbf U^{n-3})\right).
\end{equation*}
Its SSP coefficient is $\mathcal C_{\mathrm{MS}}=1/3$. 
In our implementation, the three startup states ${\mathbf U}^0,{\mathbf U}^1,{\mathbf U}^2$ are computed by the third-order SSP RK method. During these steps, we apply only the bound-preserving module at each RK stage, and the COS module after the last RK stage. For all subsequent time steps ($n\ge3$), the full EPO algorithm is used: first the positivity radius, then the multi-entropy radius, and finally the optional COS radius.
\end{example}

\section{Numerical results}\label{sec:numerics}
In this section, we present a series of numerical experiments to validate the performance of the proposed EPO discontinuous Galerkin (DG) scheme for one-dimensional (1D) and two-dimensional (2D) relativistic Euler equations.
Smooth 1D and 2D problems first verify high-order accuracy near the boundary of the invariant domain. An ultra-relativistic Riemann problem then demonstrates the effect of enforcing multiple entropy constraints. Strong-shock tests with four equations of state track the fully discrete entropy histories, and the shock--vortex and jet calculations test positivity preservation in near-vacuum and ultra-relativistic regimes.
For non-periodic problems, we monitor the boundary-corrected discrete entropy obtained by adding the accumulated numerical entropy flux through the physical boundary. The first Riemann problem is repeated with one and with three entropy pairs, and all three histories are tracked separately.

\begin{figure}[!htb]
	\centering
	\begin{subfigure}[t]{.48\textwidth}
		\centering
		\includegraphics[width=1\textwidth]{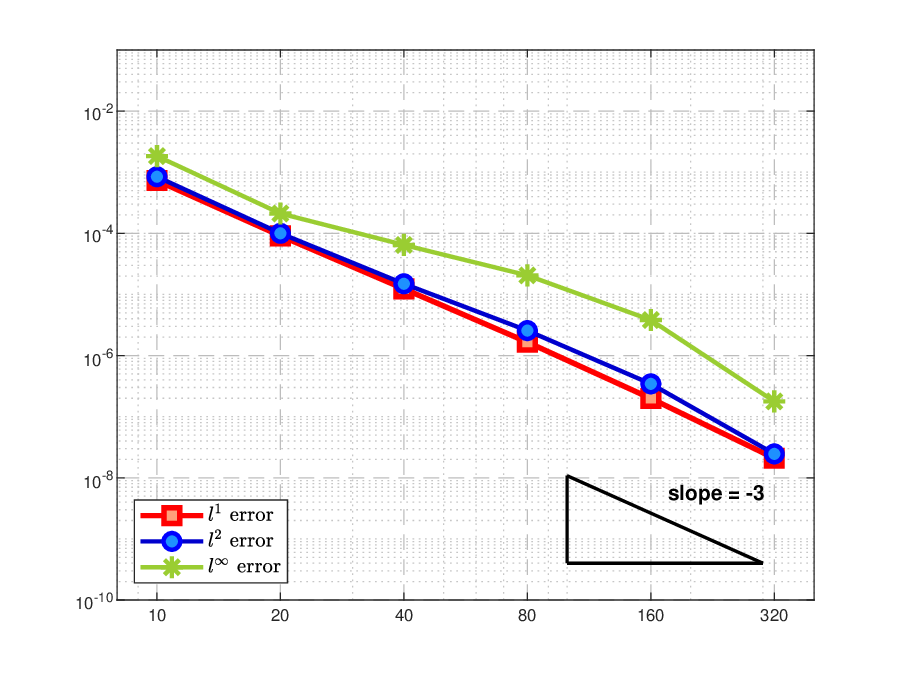}
	\end{subfigure}
    \begin{subfigure}[t]{.48\textwidth}
		\centering
		\includegraphics[width=1\textwidth]{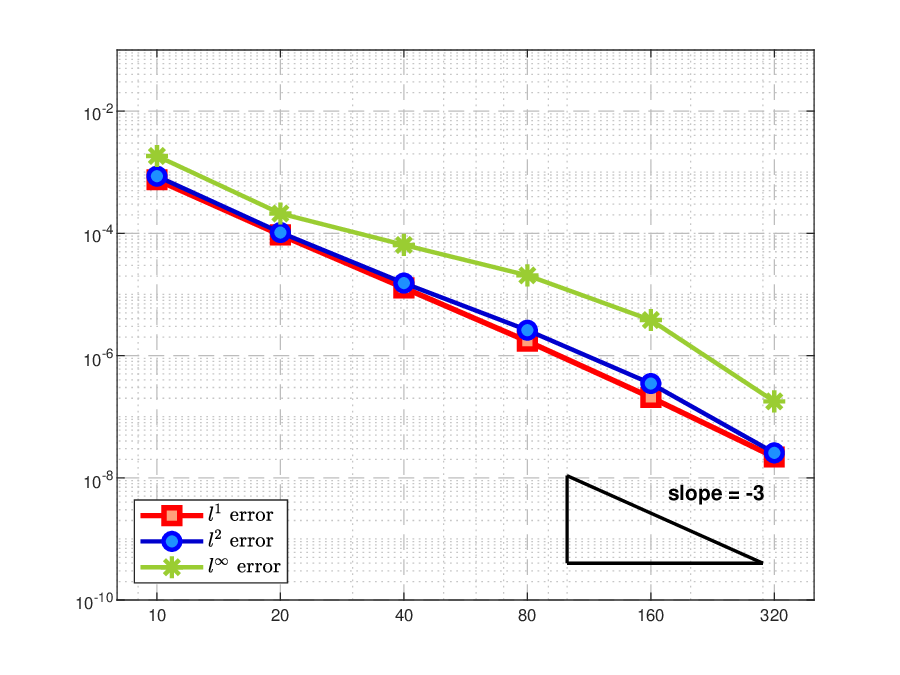}
	\end{subfigure}
	\caption{Test~\ref{Ex:1DSmooth}: The $l^1$, $l^2$, and $l^{\infty}$ errors in $\rho$ at different grid resolutions with different EOSs. Left: Ideal EOS; right: IP EOS.}\label{fig:num-smooth}
\end{figure}

All computations use a $\mathbb{P}^2$ modal DG discretization, and the Gauss--Lobatto rule in \eqref{eq:localdiscE} with $L=8$ points. The choice $L=k+1=3$ would already be exact for the polynomial degree, but the larger rule makes the local discrete entropy \eqref{eq:localdiscE} a much closer approximation of the cell entropy $|I_j|^{-1}\int_{I_j}\E(\mathbf U_h)$, so that the enforced inequality is sharper. The price is the smaller endpoint weight $\omega_1=1/\bigl(L(L-1)\bigr)=1/56$ in \eqref{eq:CFL}. The time step is chosen so that \emph{every} forward Euler step inside the third-order SSP multistep integrator of Example~\ref{ex:ssp-ms3} satisfies \eqref{eq:CFL}, i.e. $3\Delta t/\Delta x\le\omega_1$ in one dimension and
$3\Delta t\,(1/\Delta x+1/\Delta y)\le\omega_1$ in two dimensions, with a further safety factor.
The positivity thresholds are $\varepsilon_D=\varepsilon_q=10^{-13}$. Where oscillation control is applied, we use the COS module of Section~\ref{subsec:cos} with $\beta_j\equiv1$, and $C_k=0.15$ in one dimension and $C_k=0.1$ in two dimensions. The local variant of Remark~\ref{remark: local COS} uses the shock indicator with $\delta=0.05$.
Unless stated otherwise, reference solutions (for problems without closed-form exact solutions) are produced with the first-order local Lax--Friedrichs scheme on $100{,}000$ uniform cells, and are shown as solid lines, while the EPO solutions are drawn with circle markers. Four equations of state are used: ideal, RC, IP, and TM (\ref{app:eos}). Unless a different value is stated for a particular test, the ideal EOS is used with $\Gamma=5/3$.

\testcase{1D accuracy test}{Ex:1DSmooth}
We advect a smooth density wave on $[0,1]$ with periodic boundary conditions. This problem isolates the convergence behavior from any shock-capturing effect. The initial data are
\[
(\rho,v,p)(x,0) = \bigl(1 + A\sin(2\pi x),\, v_0,\, p_0\bigr),
\]
with amplitude $A = 0.99999$, constant subluminal velocity $v_0 = 0.9$, and pressure $p_0 = 1.0$. For any EOS, the exact solution is the uniform translation
\[
(\rho,v,p)(x,t) = \bigl(1 + 0.99999\sin(2\pi(x - 0.9t)),\, 0.9,\, 1\bigr).
\]
The amplitude is deliberately chosen so close to unity that the trough of the density wave is within $10^{-5}$ of vacuum, and a scheme without positivity control fails within the first few steps on coarse meshes. We run the test with the ideal and the IP equations of state.

Figure~\ref{fig:num-smooth} reports the $l^1$, $l^2$, and $l^\infty$ errors in $\rho$ against the number of cells in a log--log scale, together with a reference triangle of slope $-3$. All six curves follow the reference line under mesh refinement, so the EPO $\mathbb{P}^2$ DG scheme attains its formal third-order rate for both equations of state. The limiter does not degrade the asymptotic error, even though the solution comes within $10^{-5}$ of the boundary of the admissible set.

\begin{figure}[!htb]
	\centering
	\begin{subfigure}[t]{.48\textwidth}
		\centering
		\includegraphics[width=1\textwidth]{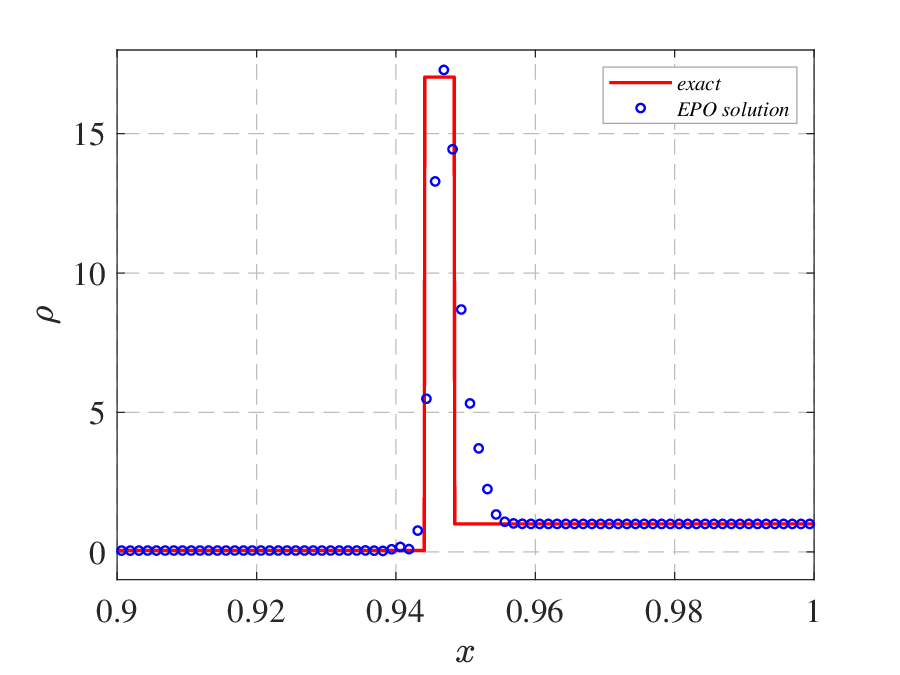}
        \caption{Close-up of $\rho$ with a single entropy pair}\label{RCSrho1EP}
	\end{subfigure}
    \begin{subfigure}[t]{.48\textwidth}
		\centering
		\includegraphics[width=1\textwidth]{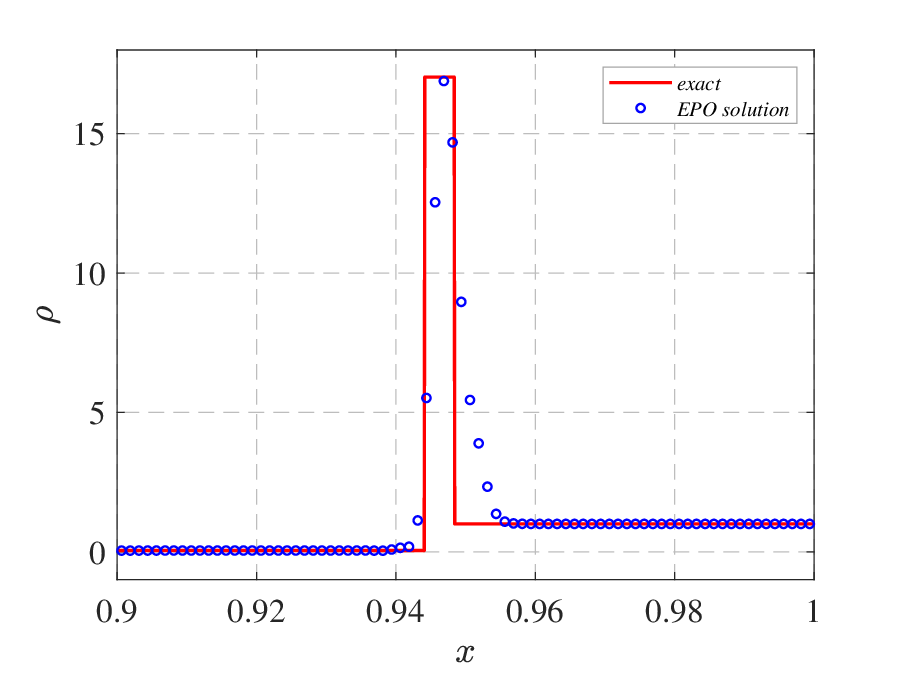}
        \caption{Close-up of $\rho$ with three entropy pairs}\label{RCSrho3EP}
	\end{subfigure}
    \begin{subfigure}[t]{.48\textwidth}
		\centering
		\includegraphics[width=1\textwidth]{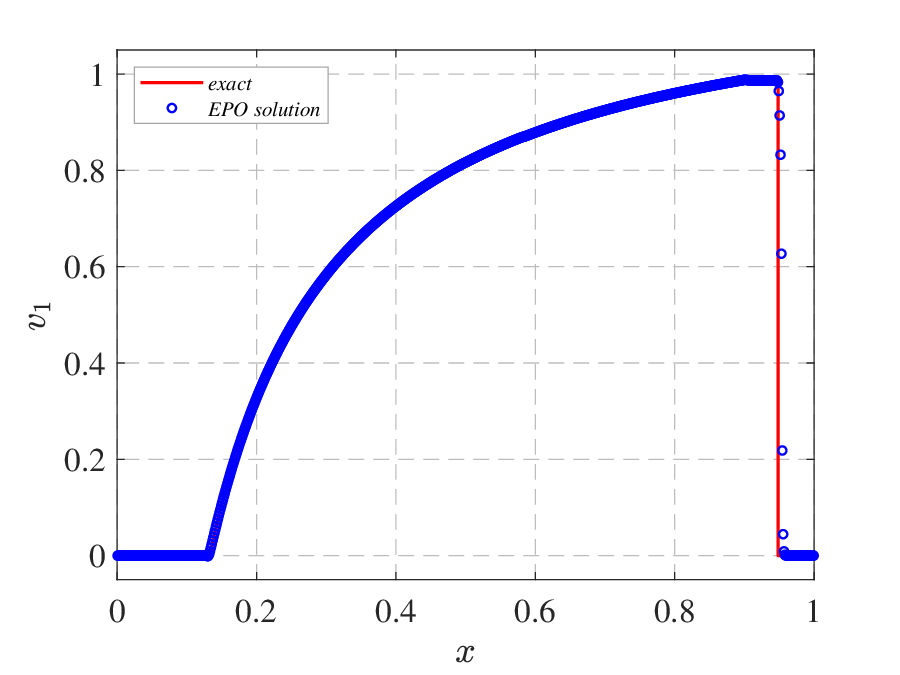}
        \caption{$v_1$}\label{RCSv13EP}
	\end{subfigure}
    \begin{subfigure}[t]{.48\textwidth}
		\centering
		\includegraphics[width=1\textwidth]{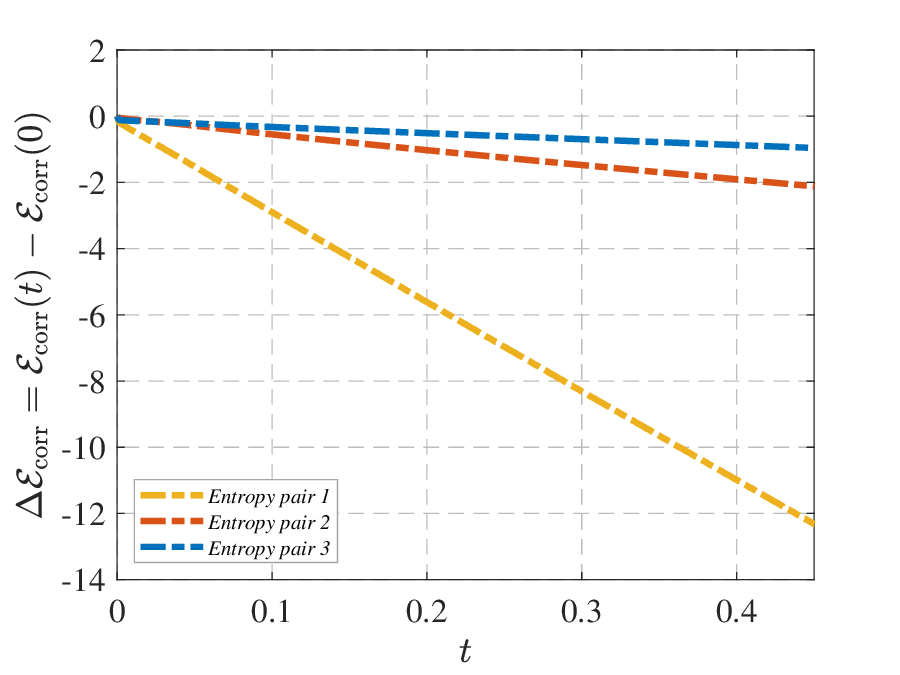}
        \caption{time evolution of discrete entropy}\label{RCSentropy}
	\end{subfigure}
	\caption{The first Riemann problem in Test~\ref{Ex:1DRP}: The EPO numerical results with multiple entropy pairs at $t = 0.45$ on $800$ uniform cells, and the time evolution of the global discrete entropy.}\label{fig:num-riemann-1}
\end{figure}

\begin{figure}[!htb]
	\centering
	\begin{subfigure}[t]{.48\textwidth}
		\centering
		\includegraphics[width=1\textwidth]{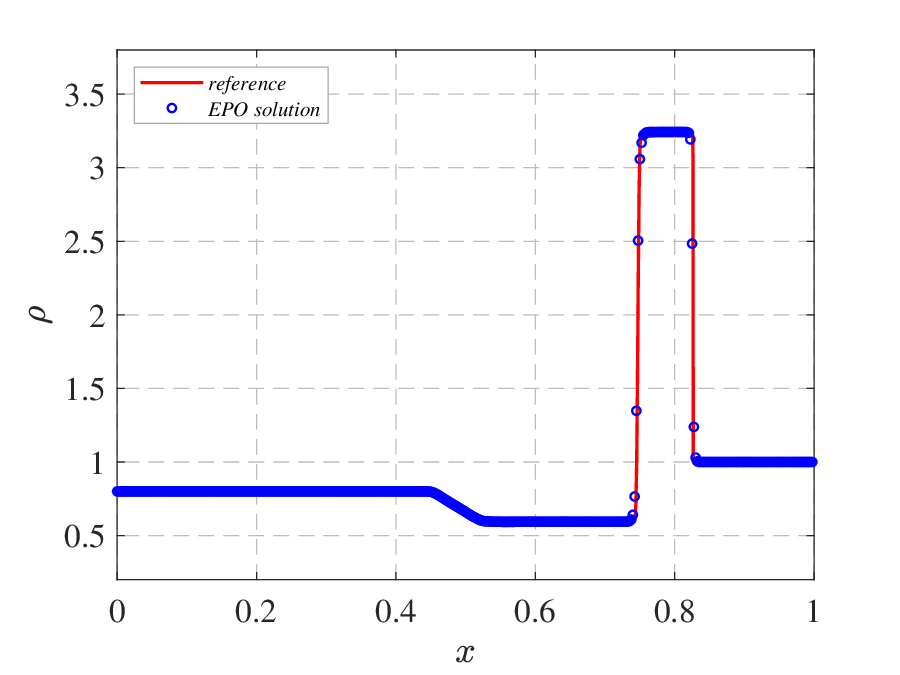}
	\end{subfigure}
    \begin{subfigure}[t]{.48\textwidth}
		\centering
		\includegraphics[width=1\textwidth]{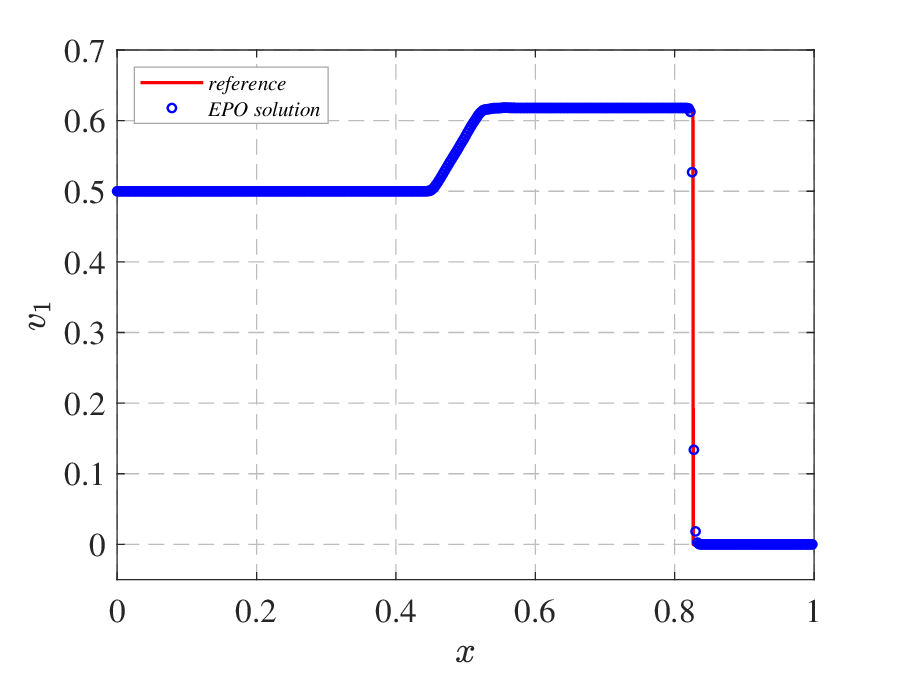}
	\end{subfigure}
    \begin{subfigure}[t]{.48\textwidth}
		\centering
		\includegraphics[width=1\textwidth]{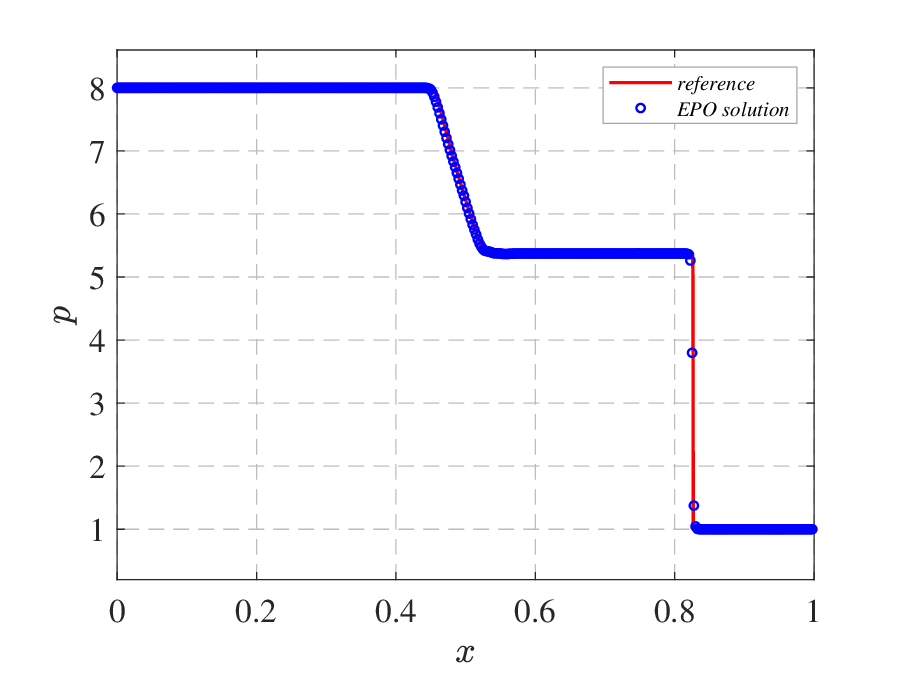}
	\end{subfigure}
        \begin{subfigure}[t]{.48\textwidth}
		\centering
		\includegraphics[width=1\textwidth]{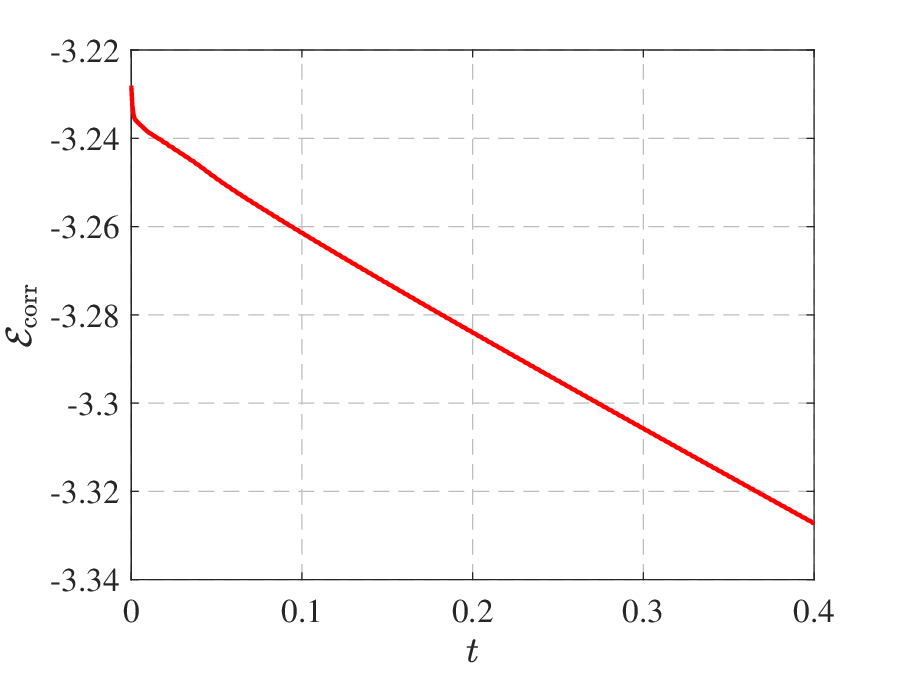}
	\end{subfigure}
	\caption{The second Riemann problem in Test~\ref{Ex:1DRP}: The EPO numerical results at $t = 0.4$ on $400$ uniform cells {with COS}, and the time evolution of the global discrete entropy.}\label{fig:num-riemann-2}
\end{figure}

\begin{figure}[!htb]
	\centering
	\begin{subfigure}[t]{.48\textwidth}
		\centering
		\includegraphics[width=1\textwidth]{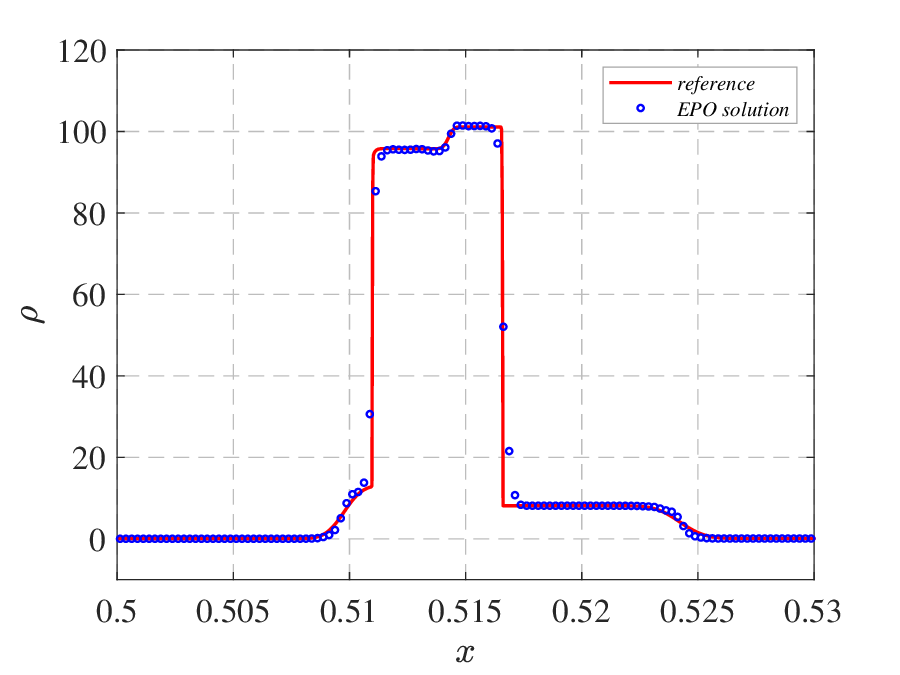}
	\end{subfigure}
    \begin{subfigure}[t]{.48\textwidth}
		\centering
		\includegraphics[width=1\textwidth]{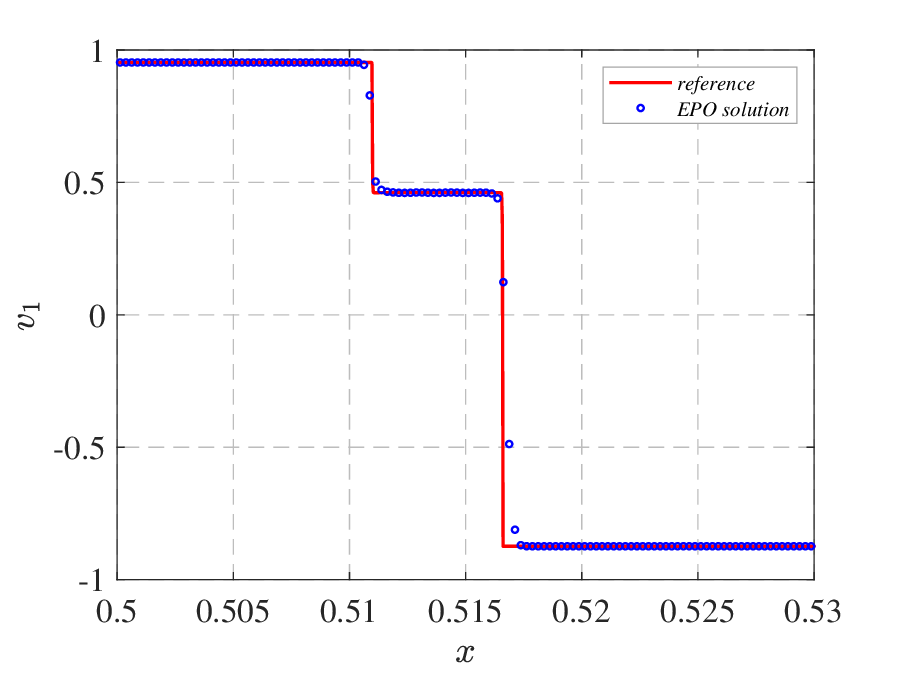}
	\end{subfigure}
    \begin{subfigure}[t]{.48\textwidth}
		\centering
		\includegraphics[width=1\textwidth]{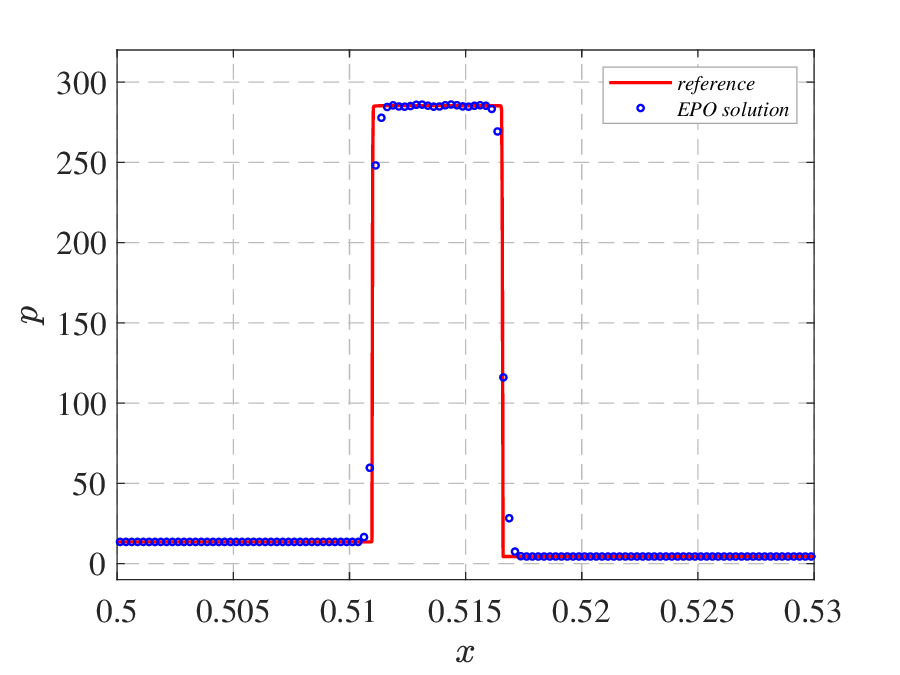}
	\end{subfigure}
    \begin{subfigure}[t]{.48\textwidth}
		\centering
		\includegraphics[width=1\textwidth]{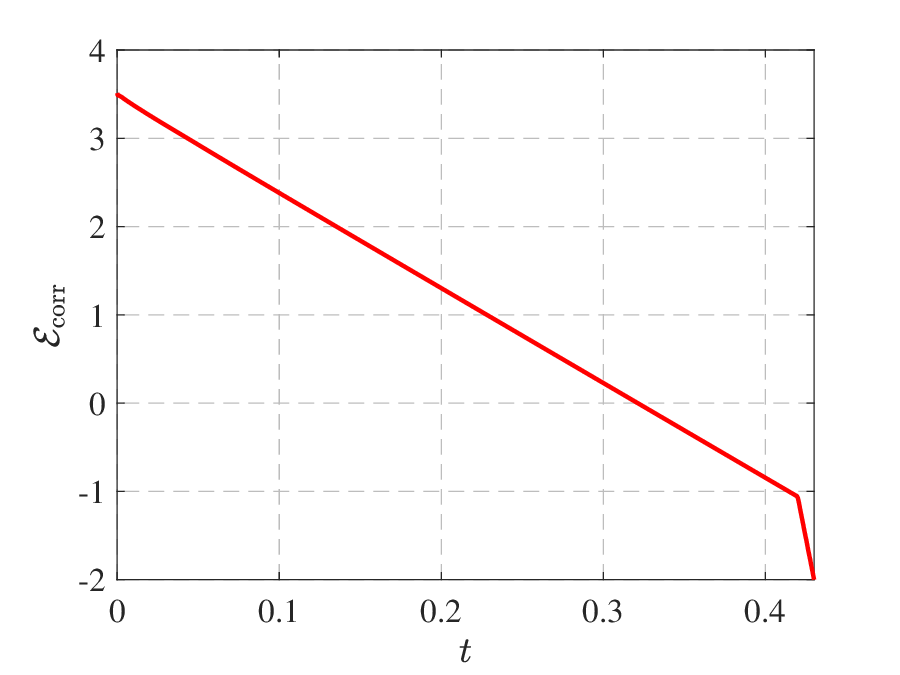}
	\end{subfigure}
	\caption{The blast problem in Test~\ref{Ex:blastwave}: The EPO numerical results at $t = 0.43$ on $4000$ uniform cells {with COS}, and the time evolution of the global discrete entropy.}\label{fig:num-BL}
\end{figure}

\testcase{1D Riemann problems}{Ex:1DRP}
Two Riemann problems are considered in this example. The first demonstrates the effect of using multiple entropy pairs, and the second exercises the full algorithm, including oscillation control with a non-ideal EOS.

\emph{(a) An ultra-relativistic Riemann problem with three entropy pairs.}
The solution consists of a strong left-going rarefaction, a right-moving contact discontinuity, and a right-moving shock. The latter two propagate at speeds extremely close to the speed of light. The initial data under the ideal EOS are
\begin{equation*}
    \textbf{V}(x,0)=
    \begin{cases}
        (1,0,10^{4})^\top, & x<0.5,\\
        (1,0,10^{-8})^\top, & x>0.5,
    \end{cases}
\end{equation*}
with outflow conditions at both boundaries. The pressure jump spans 12 orders of magnitude, and at $t=0.45$ the shell between the contact and the shock is only about $4.2\times10^{-3}$ wide. This makes the test demanding both in admissibility and resolution.

We use $800$ uniform cells and switch off the COS module, so that any oscillation control observed is attributable to the entropy constraints alone.
The computation is run twice: once with the single canonical pair \eqref{eq:canonicalpair}, and once with a family of three pairs generated by
\[
    \mathcal{H}_r(S) :=
    \begin{cases}
        \frac{\exp(r(S-S_0))-1}{r}, & r\neq0,\\[2pt]
        S-S_0, & r = 0,
    \end{cases}
\]
with reference specific entropy $S_0=0$. For the ideal EOS, the convexity condition~\eqref{eq:Hcond} requires $r < (\Gamma-1)/\Gamma$, i.e.\ $r < 2/5$ for $\Gamma = 5/3$. The values $r\in\{0,\,0.2,\,0.39\}$ are therefore permissible, and $r=0$ reproduces the canonical pair.

The results are shown in Figure~\ref{fig:num-riemann-1}. The closeups in Figures~\ref{RCSrho1EP} and \ref{RCSrho3EP} show the rest-mass density on $[0.9,1]$. With a single entropy pair, a mild overshoot remains near
$x\approx0.9469$, immediately behind the contact. With three pairs, the profile is monotone there and matches the exact solution without visible oscillation. Enforcing additional entropy inequalities acts as a mild oscillation-suppressing mechanism, with the marginal cost of two additional entropy evaluations per trial radius. Figure~\ref{RCSv13EP} shows the velocity over the whole domain computed with three pairs, and all three waves are captured at the correct locations with no smearing beyond the expected two to three cells at the shock.
Figure~\ref{RCSentropy} tracks the deviation of each boundary-corrected discrete total entropy from its initial value. All three curves leave the origin downward and remain non-increasing, confirming that the fully discrete entropy inequality holds for every pair simultaneously.

\emph{(b) A Riemann problem with the TM EOS and oscillation control.} The initial condition is
\begin{equation*}
    \textbf{V}(x,0) =
    \begin{cases}
        (0.8,0.5,8)^\top, & 0\leq x<0.5, \\
        (1,0,1)^\top, & 0.5\leq x \leq 1,
    \end{cases}
\end{equation*}
with outflow boundary conditions. The exact solution contains a left-moving rarefaction, a contact discontinuity, and a right-moving shock, all at moderate subluminal speeds.

We run the full algorithm with a single canonical entropy pair and global COS on $400$ uniform cells. Figure~\ref{fig:num-riemann-2} shows $\rho$, $v_1$ and $p$ at $t=0.4$ against the reference solution. The shock and the contact are resolved without over- or undershoots, and the rarefaction fan is clean. The entropy and COS radii together do not cause over-limiting in the smooth part of the fan. The time history of the discrete total entropy, shown in Figure~\ref{fig:num-riemann-2}, decays monotonically throughout the simulation, in agreement with the fully discrete estimate of Section~\ref{sec:epo}.

\begin{figure}[!htb]
	\centering
	\begin{subfigure}[t]{.48\textwidth}
		\centering
		\includegraphics[width=1\textwidth]{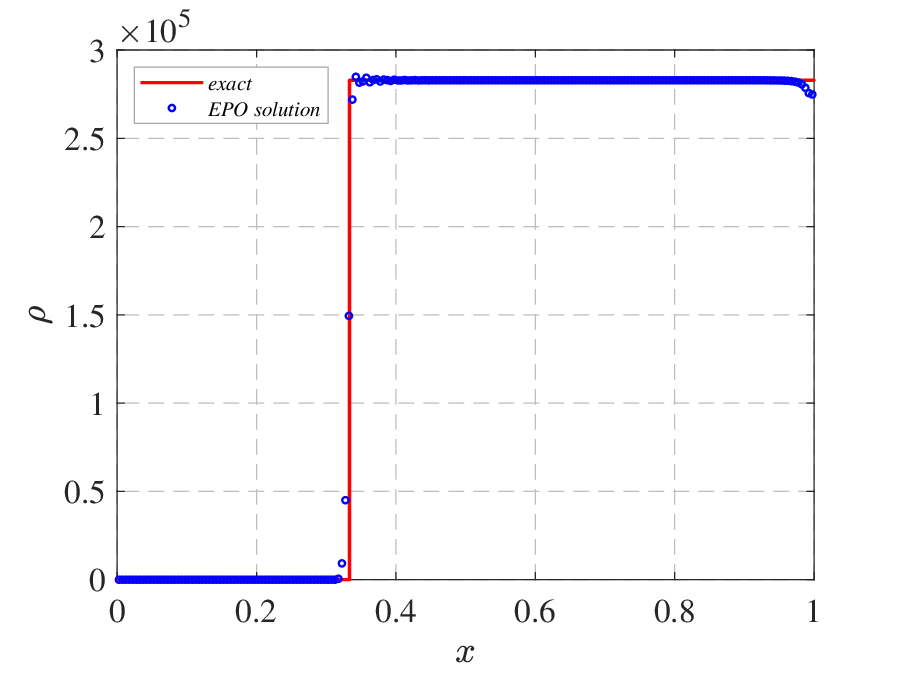}
	\end{subfigure}
    \begin{subfigure}[t]{.48\textwidth}
		\centering
		\includegraphics[width=1\textwidth]{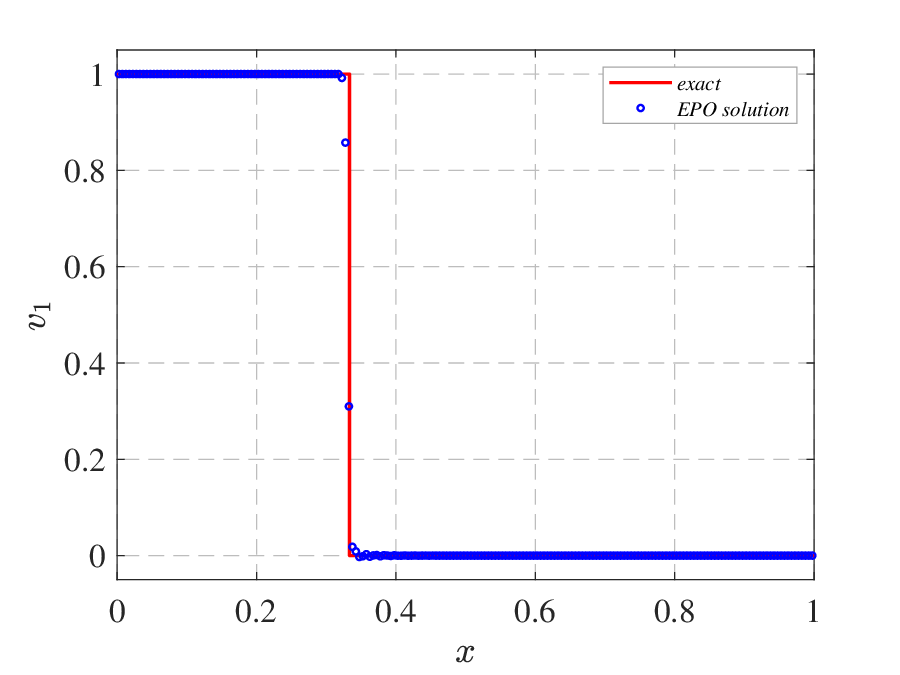}
	\end{subfigure}
    \begin{subfigure}[t]{.48\textwidth}
		\centering
		\includegraphics[width=1\textwidth]{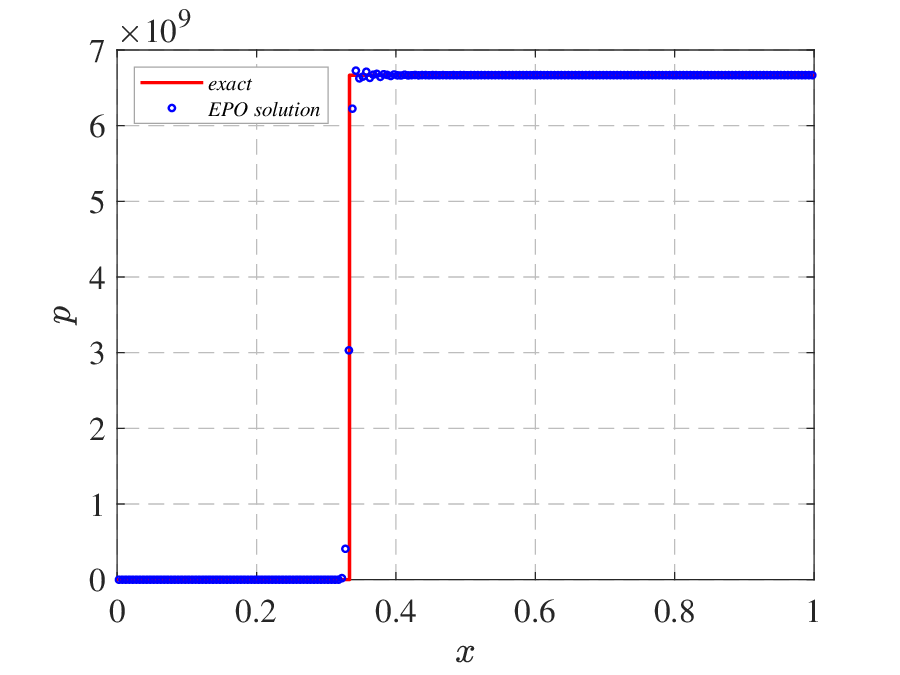}
	\end{subfigure}
    \begin{subfigure}[t]{.48\textwidth}
		\centering
		\includegraphics[width=1\textwidth]{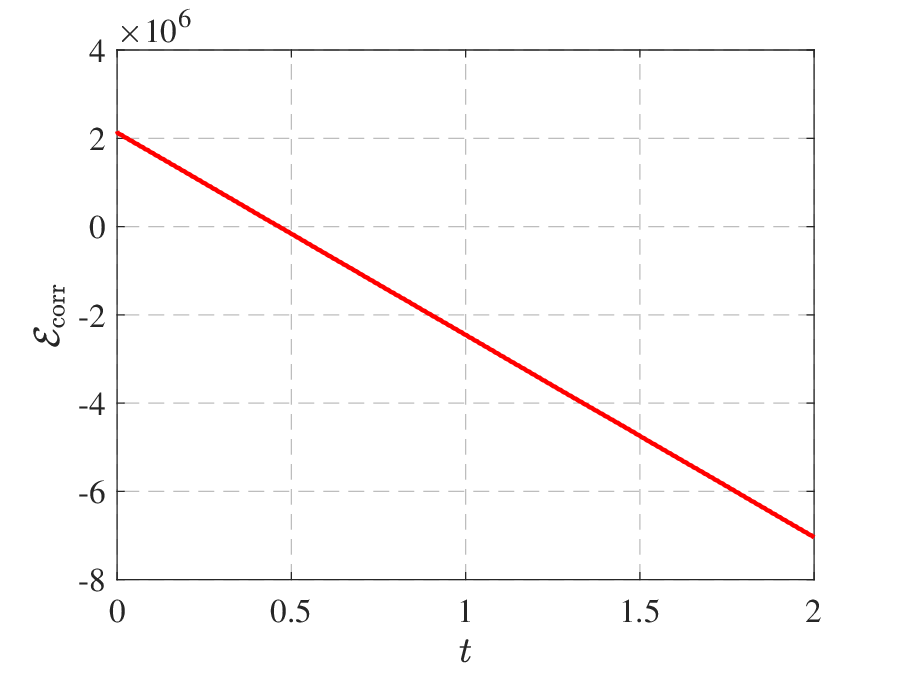}
	\end{subfigure}
	\caption{The shock-heating problem in Test~\ref{Ex:shockheating}: The EPO numerical results at $t = 2$ on $200$ uniform cells {with COS}, and the time evolution of the global discrete entropy.}\label{fig:num-SH}
\end{figure}

\testcase{Blast-wave interaction}{Ex:blastwave}
This benchmark examines the interaction between two strong relativistic blast waves, which is a stringent test for high-resolution schemes due to the formation of narrow, complex wave structures. The initial data on $[0,1]$ consist of three constant states,
\begin{equation*}
    \textbf{V}(x,0) =
    \begin{cases}
        (1,0,10^3)^\top,    & 0   \le x < 0.1,\\
        (1,0,10^{-2})^\top, & 0.1 \le x < 0.9,\\
        (1,0,10^2)^\top,    & 0.9 \le x \le 1,
    \end{cases}
\end{equation*}
with the TM EOS and outflow boundary conditions at both ends. At $t=0.43$, the solution contains two shock waves and three contact discontinuities, all compressed into the subinterval $[0.5,0.53]$. A reference solution is computed with the first-order local Lax--Friedrichs method on $200{,}000$ cells.

Figure~\ref{fig:num-BL} compares the EPO solution on $4000$ uniform cells with the reference profiles of $\rho$, $v_1$, and $p$. The density peak in the shell bounded by the two shocks, which is the most delicate feature of this test, is reproduced at the correct height and location, and no spurious extrema appear at any of the five discontinuities. The accompanying entropy history is monotone throughout the simulation, so the presence of several interacting shocks, the situation in which entropy production is largest, does not compromise the discrete entropy inequality.

\testcase{Shock heating}{Ex:shockheating}
The last 1D test models a cold, ultra-relativistic gas stream impinging on a rigid wall. This configuration produces extreme compression and a strong reflected shock. The initial data on $[0,1]$ are
\begin{equation*}
    \textbf{V}(x,0) = (1,\, 1-10^{-10},\, 10^{-4}/3)^\top,\quad 0<x<1,
\end{equation*}
with the ideal EOS, a reflective boundary at $x=1$, and an inflow boundary at $x=0$. The gas hits the wall at essentially the speed of light, and converts its kinetic energy into internal energy. A left-propagating shock forms, behind which the fluid is at rest with specific internal energy $W_0-1$, where $W_0=(1-v_0^2)^{-1/2}$ is the Lorentz factor of the incoming gas. The exact post-shock state is known analytically.

Figure~\ref{fig:num-SH} shows $\rho$, $v_1$ and $p$ at $t=2$ on $200$ uniform cells, together with the exact solution. The reflected shock is captured within two cells, and the post-shock plateau matches the analytic value. The residual wall-heating error in the pressure is confined to the two cells adjacent to the reflecting boundary. With an initial pressure of order $10^{-4}$ and an incoming Lorentz factor of order $10^{5}$, the compressed state sits very close to the boundary of $\Gset$, but no inadmissible nodal value occurs during the simulation. The entropy history in the same figure decreases at every step.

\begin{table}[!thb]
    \renewcommand{\arraystretch}{1.5}
    \centering
    \belowrulesep=0pt
    \aboverulesep=0pt
    \caption{Test~\ref{Ex:2DSmooth}: Numerical errors and convergence rates in $\rho$ at different grid resolutions.}
    \label{table:2DSmooth}
    \setlength{\tabcolsep}{2mm}{
        \begin{tabular}{c|c|ccccccc}
            \toprule[1.5pt]
            \multirow{2}{*}{EOS } &
            \multirow{2}{*}{$N_x\times N_y$} &
            \multicolumn{2}{c}{$l^1$ norm} &
            \multicolumn{2}{c}{$l^2$ norm} &
            \multicolumn{2}{c}{$l^{\infty}$ norm}\\
            \cmidrule(r){3-4} \cmidrule(r){5-6} \cmidrule(r){7-8}
            & & error & order &  error & order &  error & order\\
            \midrule[1.5pt]
            \multirow{6}{*}{Ideal} & $10\times10$ & 4.49e-03 & -- &  5.98e-03 & -- &  2.71e-02 & --    \\
            & $20\times20$ & 3.65e-04 & 3.62 & 4.83e-04 & 3.63 &  1.99e-03 & 3.77 \\
            & $40\times40$ & 4.13e-05 & 3.14 & 5.53e-05 & 3.13 &  2.06e-04 & 3.28 \\
            & $80\times80$ & 5.26e-06 & 2.98 & 7.01e-06 & 2.98 &  2.56e-05 & 3.00 \\
            & $160\times160$ & 6.47e-07 & 3.02 & 8.76e-07 & 3.00 &  3.20e-06 & 3.00 \\
            & $320\times320$ & 7.64e-08 & 3.08 & 1.05e-07 & 3.06 &  4.00e-07 & 3.00 \\
            \midrule
            \multirow{6}{*}{RC} & $10\times10$ & 4.61e-03 & -- &  6.18e-03 & -- &  2.78e-02 & --    \\
            & $20\times20$ & 4.27e-04 & 3.43 & 5.52e-04 & 3.48 &  2.21e-03 & 3.65 \\
            & $40\times40$ & 5.95e-05 & 2.84 & 8.33e-05 & 2.73 &  3.47e-04 & 2.67 \\
            & $80\times80$ & 7.41e-06 & 3.00 & 1.09e-05 & 2.93 &  6.75e-05 & 2.36 \\
            & $160\times160$ & 7.54e-07 & 3.30 & 1.03e-06 & 3.41 &  4.94e-06 & 3.77 \\
            & $320\times320$ & 8.26e-08 & 3.19 & 1.10e-07 & 3.23 &  4.00e-07 & 3.63 \\
            \bottomrule[1.5pt]
        \end{tabular}
    }
\end{table}

\begin{figure}[!htb]
	\centering
	\begin{subfigure}[t]{.48\textwidth}
		\centering
		\includegraphics[width=1\textwidth]{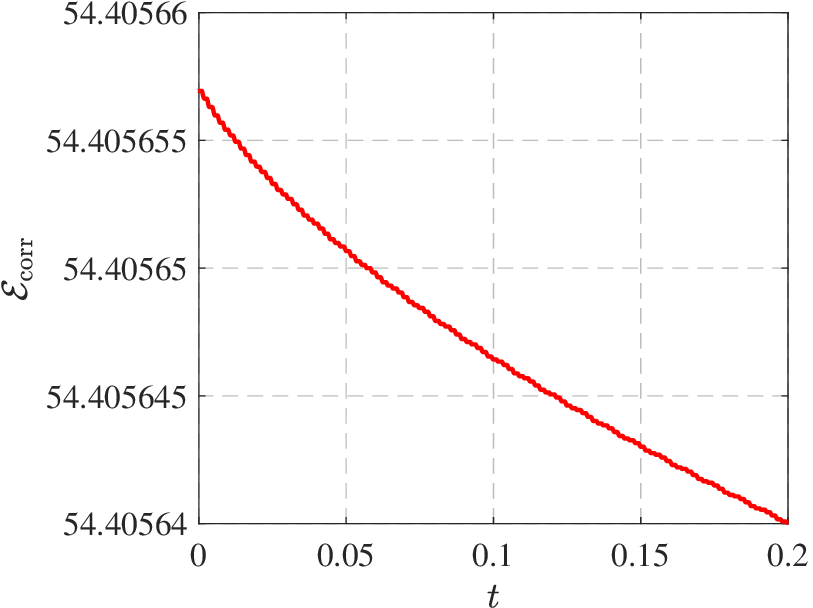}
	\end{subfigure}
    \begin{subfigure}[t]{.48\textwidth}
		\centering
		\includegraphics[width=1\textwidth]{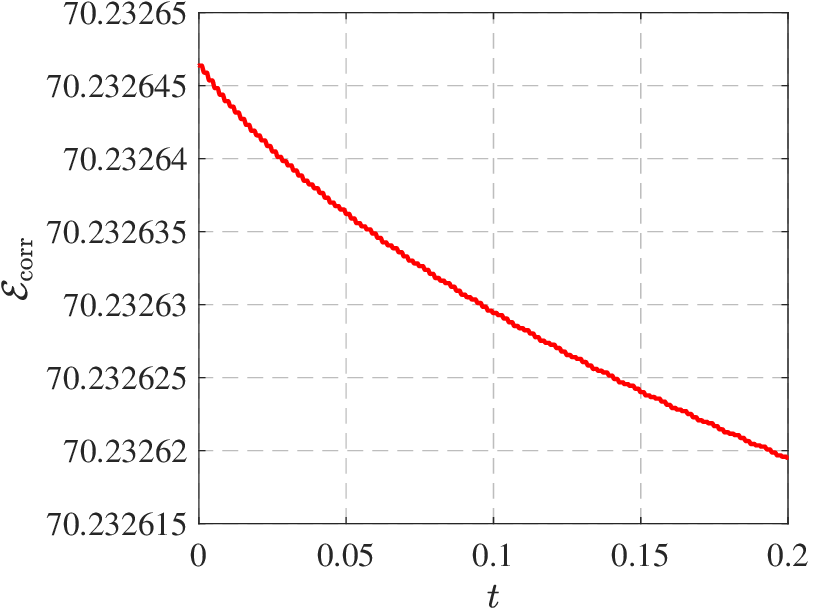}
	\end{subfigure}
	\caption{Test~\ref{Ex:2DSmooth}: The time evolution of the global discrete entropy for the resolution $N_x\times N_y = 80\times80$. Left: Ideal EOS; right: RC EOS.}\label{fig:num-2Dsmooth}
\end{figure}

\testcase{2D smooth problem}{Ex:2DSmooth}
We consider the 2D version of Test~\ref{Ex:1DSmooth}. A smooth wave propagates diagonally across the periodic unit square $[0,1]^2$, with initial data
\[
(\rho,v_x,v_y,p)(x,y,0)
= \bigl(1 + A\sin(2\pi(x+y)),\, v_x^0,\, v_y^0,\, p_0\bigr),
\]
amplitude $A = 0.99999$, velocities $(v_x^0,v_y^0) = (0.99/\sqrt{2}, 0.99/\sqrt{2})$ and pressure $p_0 = 0.01$. Both the ideal and the RC equations of state are used. The low ambient pressure, combined with the near-unit amplitude, brings this test close to both vacuum and light-speed constraints.

Table~\ref{table:2DSmooth} lists the $l^1$, $l^2$, and $l^\infty$ errors in $\rho$ on a sequence of uniformly refined Cartesian meshes with $\mathbb{P}^2$ elements, together with the observed rates. All rates approach $3$ under refinement for both equations of state, so the projection of Section~\ref{sec:epo} preserves accuracy. The mild non-monotonicity of the RC rates on the coarse meshes is a pre-asymptotic effect and is not caused by the limiter, since the errors themselves decrease monotonically.
Figure~\ref{fig:num-2Dsmooth} shows the global discrete entropy on the $80\times80$ mesh. It decreases monotonically for both equations of state, as predicted by the estimate of Section~\ref{sec:epo}. The exact solution produces no entropy, and the observed decay is purely numerical. Its magnitude measures the dissipation introduced by the $\alpha=1$ interface flux.

\begin{figure}[!htb]
	\centering
	\begin{subfigure}[t]{.48\linewidth}
		\centering
		\includegraphics[width=0.95\textwidth]{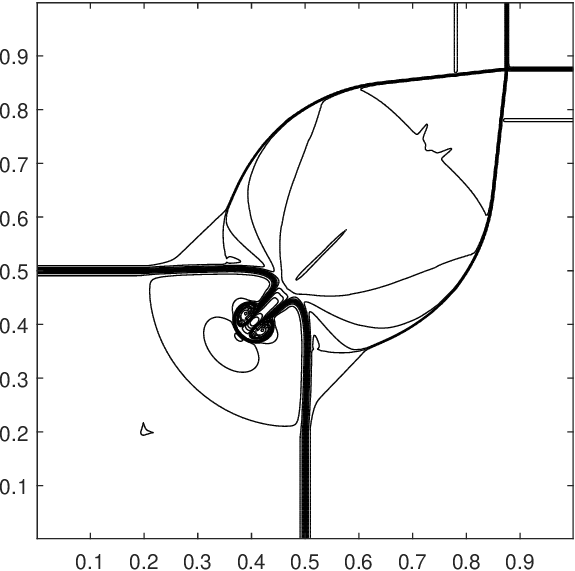}
	\end{subfigure}
	\begin{subfigure}[t]{.48\linewidth}
		\centering
		\includegraphics[width=0.95\textwidth]{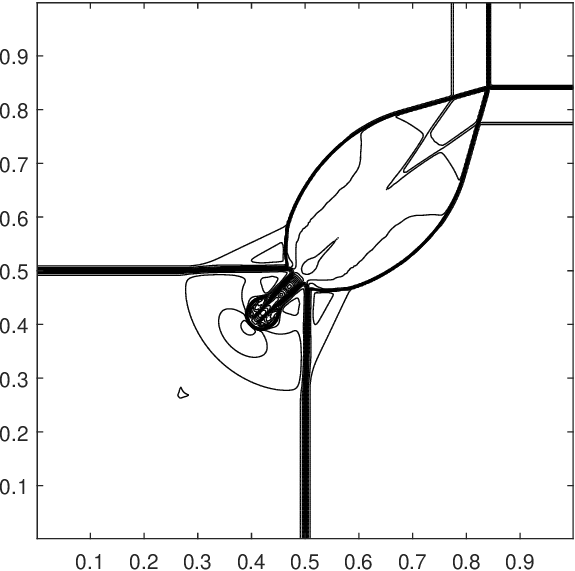}
	\end{subfigure}
    \begin{subfigure}[t]{.48\linewidth}
		\centering
		\includegraphics[width=0.95\textwidth]{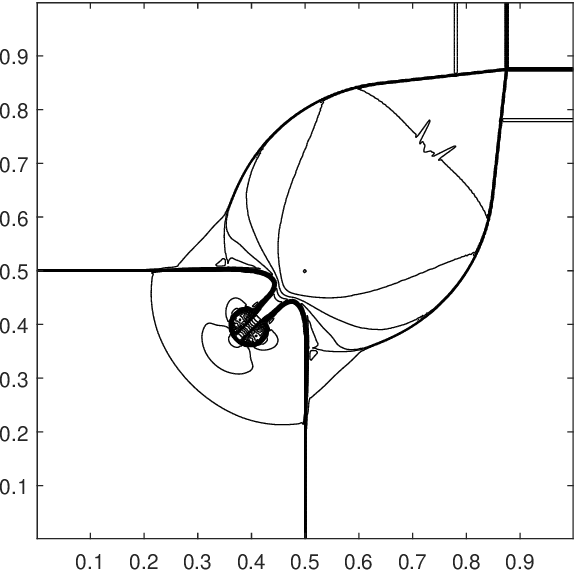}
	\end{subfigure}
	\begin{subfigure}[t]{.48\linewidth}
		\centering
		\includegraphics[width=0.95\textwidth]{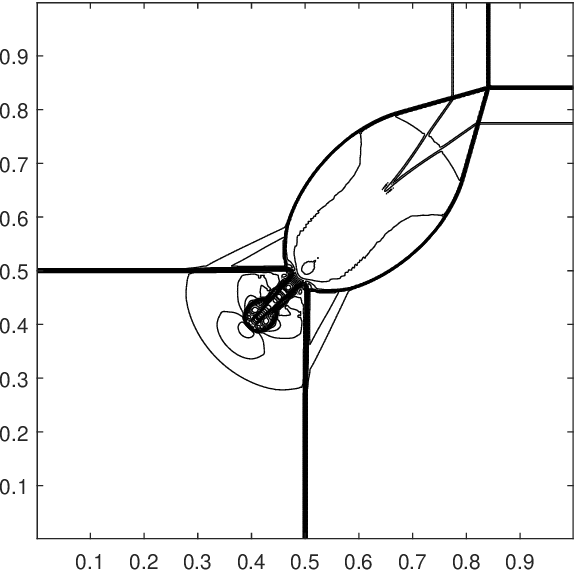}
	\end{subfigure}
	\caption{The first problem of Test~\ref{Ex:2DRP}: The contours of $\ln\rho$ on $400\times400$ uniform cells at $t = 0.4$ with (top) global COS or (bottom) local COS. Twenty-five equally spaced contour lines from -3.3 to -0.4 are displayed. Left: Ideal EOS; right: RC EOS.}\label{fig:num-2Driemann-T1}
\end{figure}
\begin{figure}[!htb]
	\centering
	\begin{subfigure}[t]{.48\linewidth}
		\centering
		\includegraphics[width=0.95\textwidth]{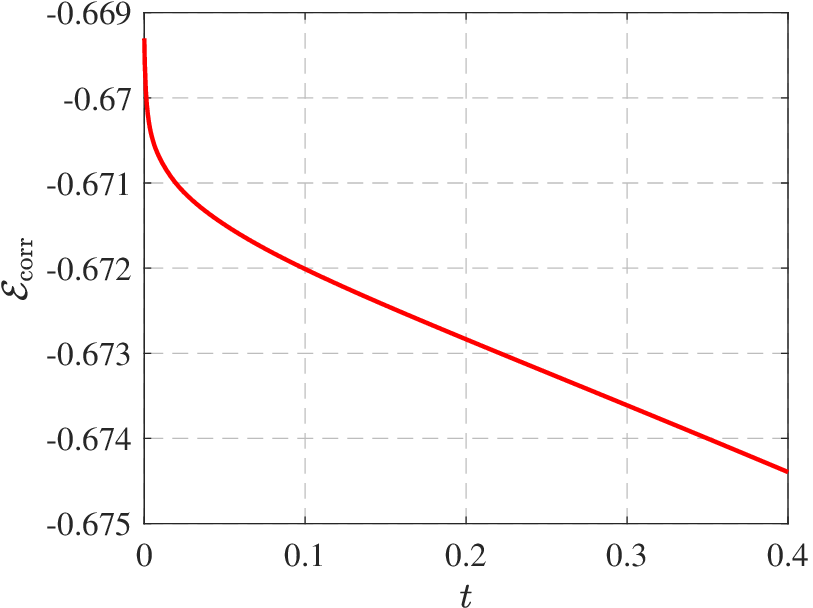}
	\end{subfigure}
	\begin{subfigure}[t]{.48\linewidth}
		\centering
		\includegraphics[width=0.95\textwidth]{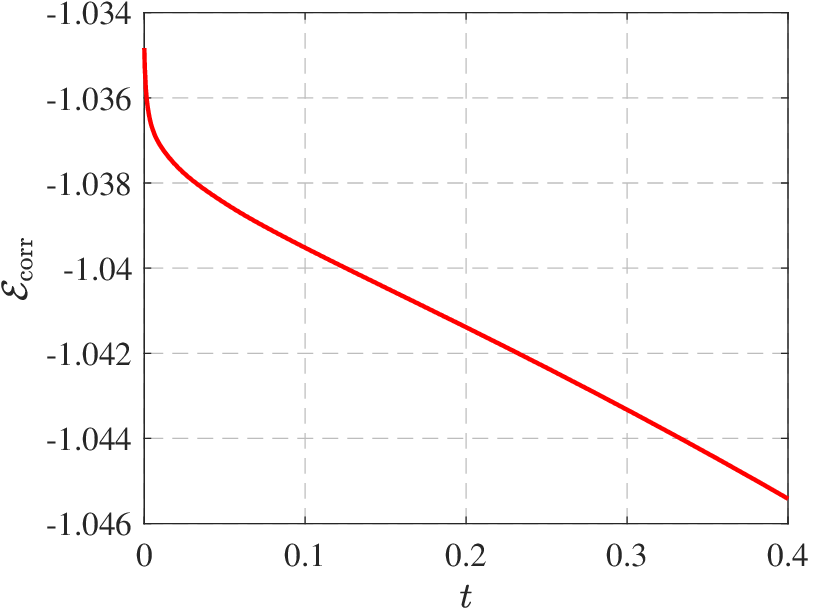}
	\end{subfigure}
	\caption{The first problem of Test~\ref{Ex:2DRP}: The time evolution of the global discrete entropy. Left: Ideal EOS; right: RC EOS.}\label{fig:num-2Driemann-T1-entropy}
\end{figure}

\begin{figure}[!htb]
	\centering
	\begin{subfigure}[t]{.48\linewidth}
		\centering
		\includegraphics[width=0.95\textwidth]{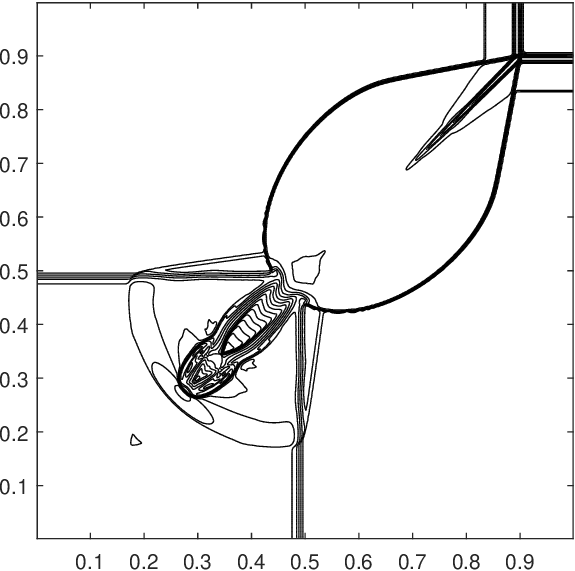}
	\end{subfigure}
	\begin{subfigure}[t]{.48\linewidth}
		\centering
		\includegraphics[width=0.95\textwidth]{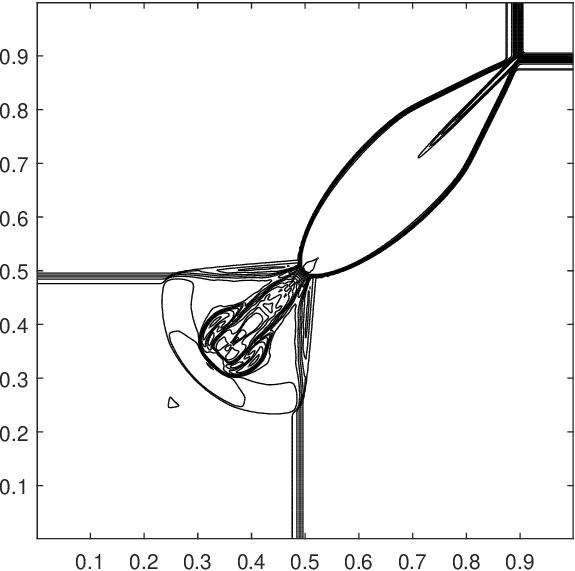}
	\end{subfigure}
	\caption{The second problem of Test~\ref{Ex:2DRP}: The contours of $\ln\rho$ on $400\times400$ uniform cells at $t = 0.4$ {with COS}. Twenty-five equally spaced contour lines from -6.0 to 1.9 are displayed. Left: Ideal EOS; right: TM EOS.}\label{fig:num-2Driemann-T4}
\end{figure}
\begin{figure}[!htb]
	\centering
	\begin{subfigure}[t]{.48\linewidth}
		\centering
		\includegraphics[width=0.95\textwidth]{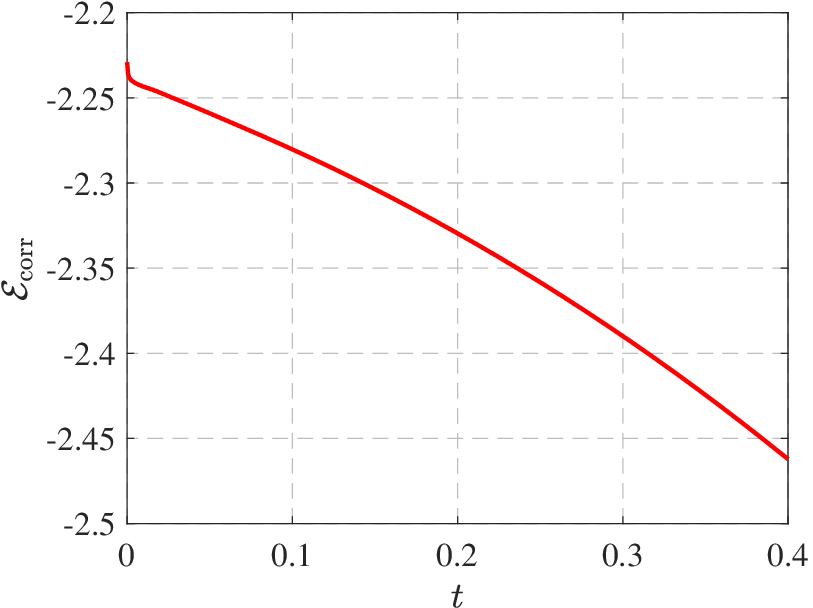}
	\end{subfigure}
	\begin{subfigure}[t]{.48\linewidth}
		\centering
		\includegraphics[width=0.95\textwidth]{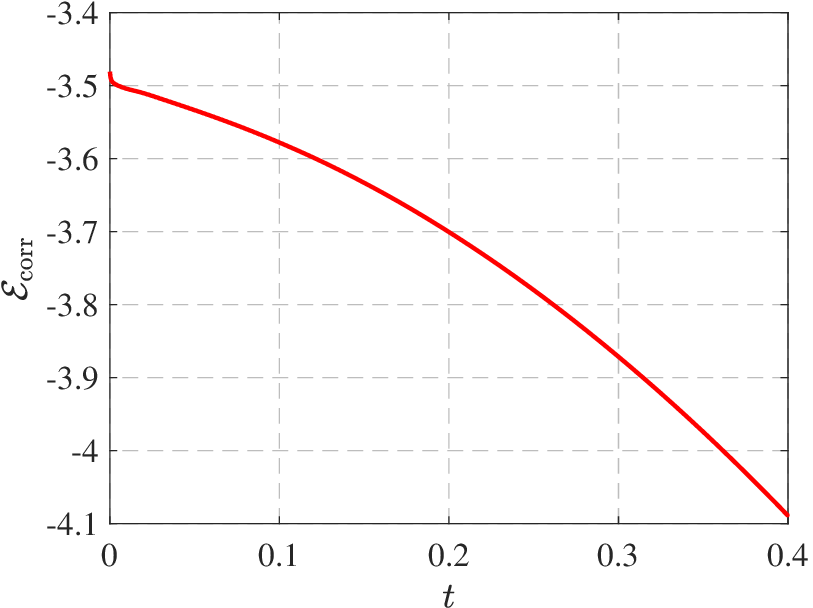}
	\end{subfigure}
	\caption{The second problem of Test~\ref{Ex:2DRP}: The time evolution of the global discrete entropy. Left: Ideal EOS; right: TM EOS.}\label{fig:num-2Driemann-T4-entropy}
\end{figure}

\begin{table}[!thb]
    \renewcommand{\arraystretch}{1.5}
    \centering
    \belowrulesep=0pt
    \aboverulesep=0pt
    \caption{Initial data of the four 2D RPs in Test~\ref{Ex:2DRP}.}
    \label{table:2DRP}
    \scalebox{0.62}{
    \setlength{\tabcolsep}{2mm}{
        \begin{tabular}{c|c|cccc}
            \toprule[1.5pt]
            {RP } &
            {domain} &
            {$\rho$} &
            {$v_1$} &
            {$v_2$} &
            {$p$}\\
            \midrule[1.5pt]
            \multirow{4}{*}{RP I}
            & $x>0.5,\,y>0.5$ & 0.035145216124503 & 0 &  0 & 0.162931056509027  \\
            \cmidrule(r){2-6}
            & $x<0.5,\,y>0.5$ & 0.1 & 0.7 & 0 & 1  \\
            \cmidrule(r){2-6}
            & $x<0.5,\,y<0.5$ & 0.5 & 0 & 0 & 1  \\
            \cmidrule(r){2-6}
            & $x>0.5,\,y<0.5$ & 0.1 & 0 & 0.7 & 1  \\
            \midrule
            \multirow{4}{*}{RP II}
            & $x>0.5,\,y>0.5$ & 0.1 & 0 &  0 & 0.01  \\
            \cmidrule(r){2-6}
            & $x<0.5,\,y>0.5$ & 0.1 & 0.99 & 0 & 1 \\
            \cmidrule(r){2-6}
            & $x<0.5,\,y<0.5$ & 0.5 & 0 & 0 & 1 \\
            \cmidrule(r){2-6}
            & $x>0.5,\,y<0.5$ & 0.1 & 0 & 0.99 & 1 \\
            \midrule
            \multirow{4}{*}{RP III}
            & $x>0.5,\,y>0.5$ & 0.5 & 0.5 &  -0.5 & 5  \\
            \cmidrule(r){2-6}
            & $x<0.5,\,y>0.5$ & 1 & 0.5 & 0.5 & 5  \\
            \cmidrule(r){2-6}
            & $x<0.5,\,y<0.5$ & 3 & -0.5 & 0.5 & 5  \\
            \cmidrule(r){2-6}
            & $x>0.5,\,y<0.5$ & 1.5 & -0.5 & -0.5 & 5  \\
            \midrule
            \multirow{4}{*}{RP IV}
            & $x>0.5,\,y>0.5$ & 0.1 & 0 &  0 & 20  \\
            \cmidrule(r){2-6}
            & $x<0.5,\,y>0.5$ & 0.00414329639576 & 0.9946418833556542 & 0 & 0.05  \\
            \cmidrule(r){2-6}
            & $x<0.5,\,y<0.5$ & 0.01 & 0 & 0 & 0.05  \\
            \cmidrule(r){2-6}
            & $x>0.5,\,y<0.5$ & 0.00414329639576 & 0 & 0.9946418833556542 & 0.05  \\
            \bottomrule[1.5pt]
        \end{tabular}
    }
    }
\end{table}

\testcase{2D four-state Riemann problems}{Ex:2DRP}
This set of benchmarks consists of four classical 2D Riemann problems. Each starts from four constant states in the four quadrants of the unit square $[0,1]^2$, as listed in Table~\ref{table:2DRP}. We use a uniform $400\times400$ mesh with outflow boundary conditions, and test each configuration with the ideal EOS and one additional EOS: RC for RP~I and RP~IV, TM for RP~II, and IP for RP~III. The four problems exercise different dynamics: RP~I produces a mushroom-shaped roll-up from the interaction of two contacts and two shocks, RP~II is its ultra-relativistic variant and develops a thin diagonal jet, RP~III generates a spiral fed by four contact discontinuities, and RP~IV combines velocities close to the speed of light with a strong pressure contrast.

\begin{figure}[!htb]
	\centering
	\begin{subfigure}[t]{.48\linewidth}
		\centering
		\includegraphics[width=0.95\textwidth]{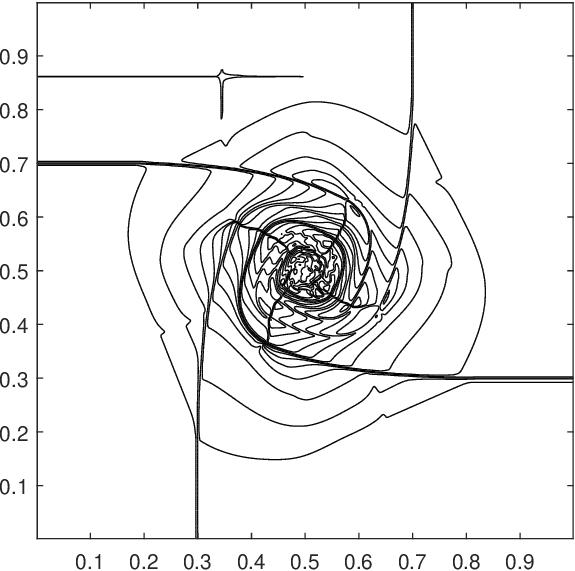}
	\end{subfigure}
	\begin{subfigure}[t]{.48\linewidth}
		\centering
		\includegraphics[width=0.95\textwidth]{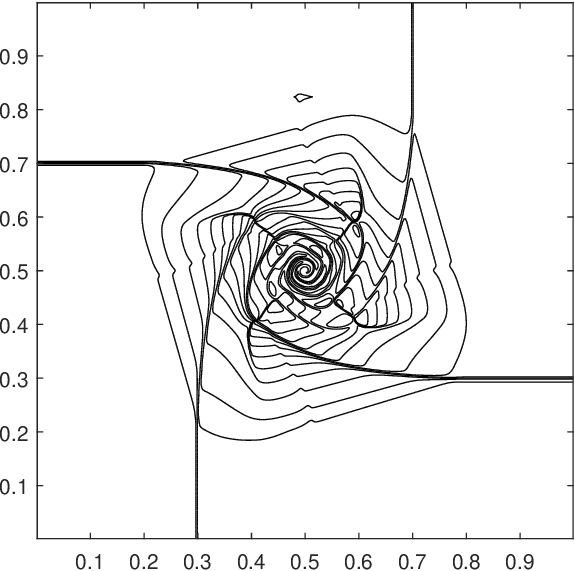}
	\end{subfigure}
	\caption{The third problem of Test~\ref{Ex:2DRP}: The contours of $\ln\rho$ on $400\times400$ uniform cells at $t = 0.4$ {with COS}. Twenty-five equally spaced contour lines from -6.0 to 1.3 are displayed. Left: Ideal EOS; right: IP EOS.}\label{fig:num-2Driemann-T2}
\end{figure}

\begin{figure}[!htb]
	\centering
	\begin{subfigure}[t]{.48\linewidth}
		\centering
		\includegraphics[width=0.95\textwidth]{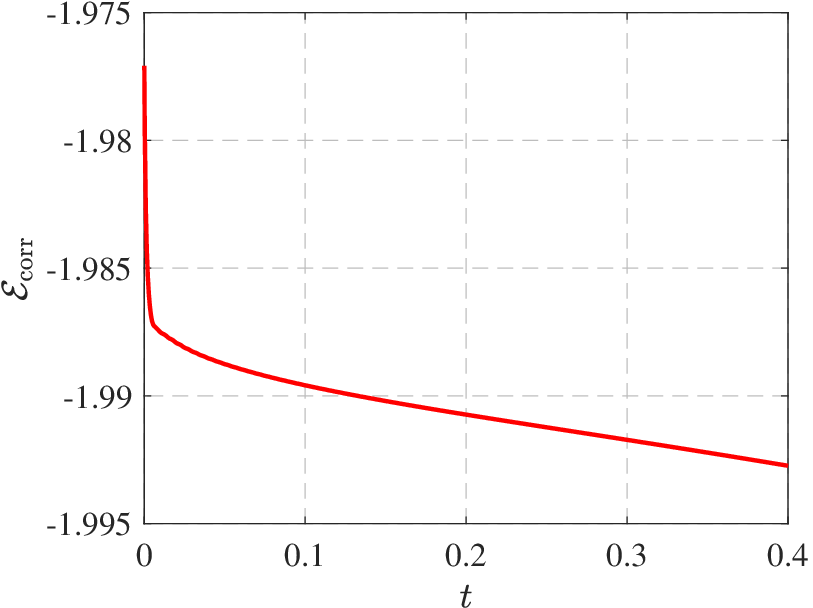}
	\end{subfigure}
	\begin{subfigure}[t]{.48\linewidth}
		\centering
		\includegraphics[width=0.95\textwidth]{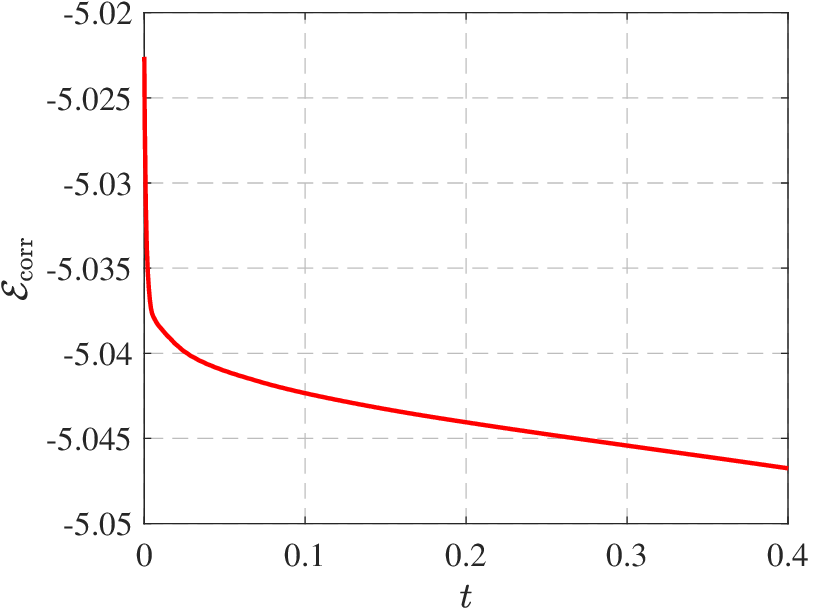}
	\end{subfigure}
	\caption{The third problem of Test~\ref{Ex:2DRP}: The time evolution of the global discrete entropy. Left: Ideal EOS; right: IP EOS.}\label{fig:num-2Driemann-T2-entropy}
\end{figure}

\begin{figure}[!htb]
	\centering
	\begin{subfigure}[t]{.48\linewidth}
		\centering
		\includegraphics[width=0.95\textwidth]{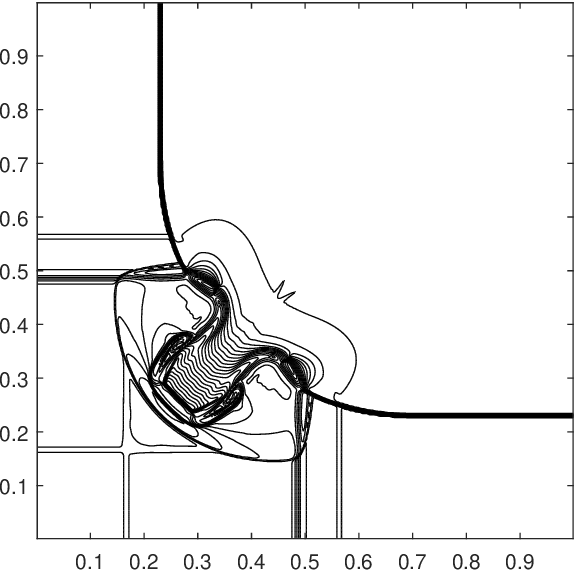}
	\end{subfigure}
	\begin{subfigure}[t]{.48\linewidth}
		\centering
		\includegraphics[width=0.95\textwidth]{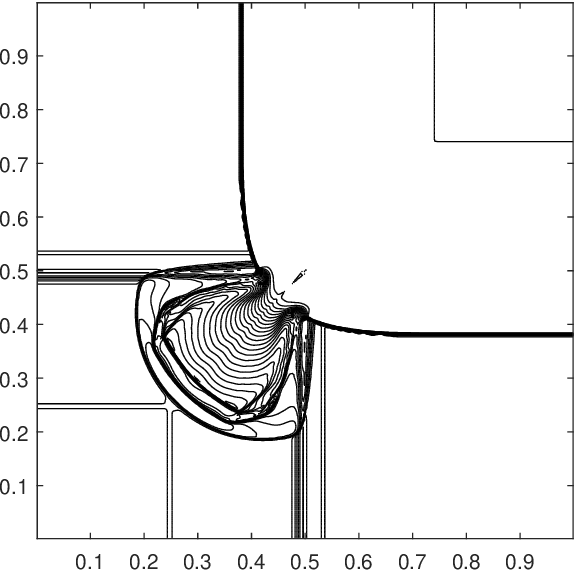}
	\end{subfigure}
	\caption{The fourth problem of Test~\ref{Ex:2DRP}: The contours of $\ln\rho$ on $400\times400$ uniform cells at $t = 0.4$. Twenty-five equally spaced contour lines from -8.0 to -2.3 are displayed. Left: Ideal EOS; right: RC EOS.}\label{fig:num-2Driemann-T5}
\end{figure}

\begin{figure}[!htb]
	\centering
	\begin{subfigure}[t]{.48\linewidth}
		\centering
		\includegraphics[width=0.95\textwidth]{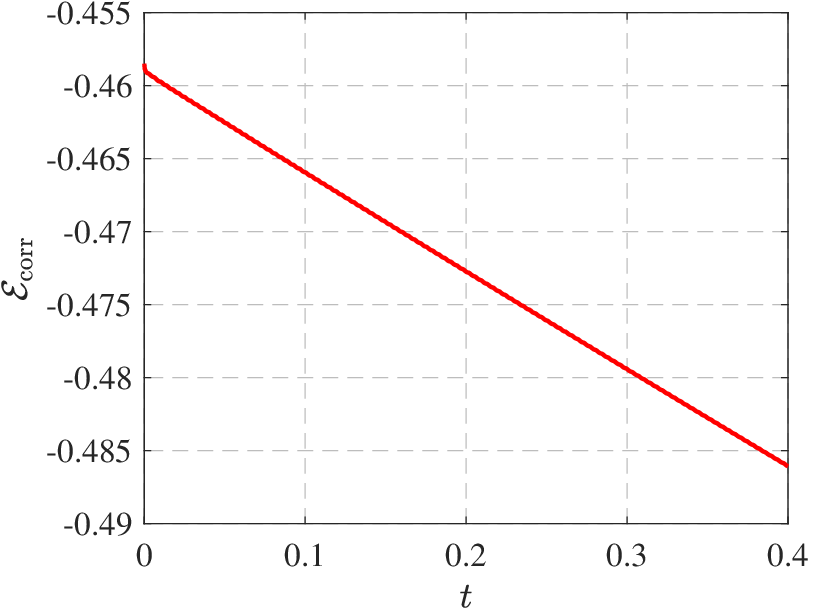}
	\end{subfigure}
	\begin{subfigure}[t]{.48\linewidth}
		\centering
		\includegraphics[width=0.95\textwidth]{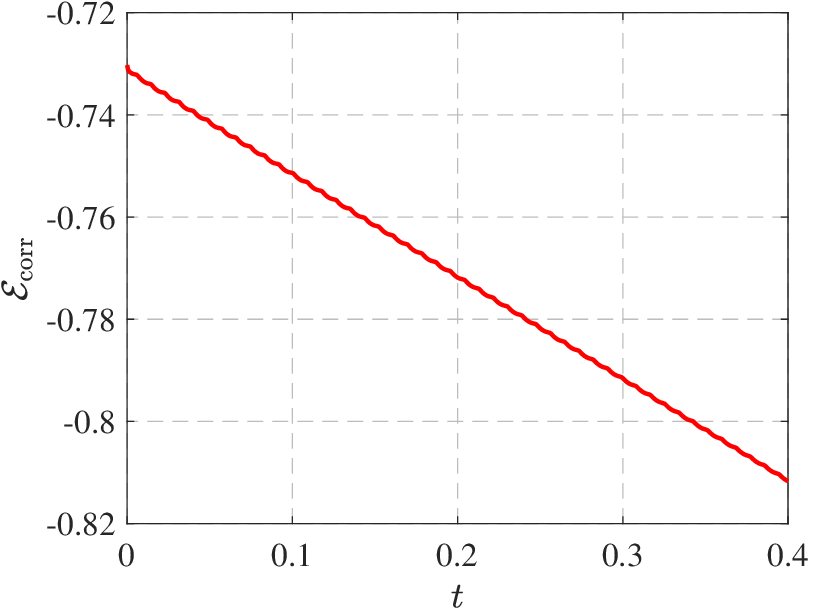}
	\end{subfigure}
	\caption{The fourth problem of Test~\ref{Ex:2DRP}: The time evolution of the global discrete entropy. Left: Ideal EOS; right: RC EOS.}\label{fig:num-2Driemann-T5-entropy}
\end{figure}

Figures~\ref{fig:num-2Driemann-T1}, \ref{fig:num-2Driemann-T4},
\ref{fig:num-2Driemann-T2} and \ref{fig:num-2Driemann-T5} show the $\ln\rho$ contours at $t=0.4$ for RP~I--IV, and
Figures~\ref{fig:num-2Driemann-T1-entropy}, \ref{fig:num-2Driemann-T4-entropy}, \ref{fig:num-2Driemann-T2-entropy} and \ref{fig:num-2Driemann-T5-entropy} display the corresponding entropy histories. For every configuration and every EOS, the computed state remains admissible at all quadrature nodes throughout the simulation. This holds even in RP~II and RP~IV, where the maximal fluid velocity is within $10^{-2}$ of the speed of light. The discrete total entropy decreases monotonically in all cases. The contour plots capture
the mushroom roll-up, the reflected curved shocks, the diagonal jet, and the spiral without visible spurious oscillations.

RP~I also illustrates a trade-off. With the global COS module, the two contact discontinuities that bound the mushroom stem are visibly smeared. This is because the indicator \eqref{eq:COSsigma} damps every cell where the solution is non-smooth, including contacts. Replacing it with the local variant of Remark~\ref{remark: local COS}, which restricts the damping to interfaces flagged by the Lax entropy criterion, recovers the contacts while keeping the
shocks equally clean, as the two rows of Figure~\ref{fig:num-2Driemann-T1} show. Since the oscillation radius is computed independently of the entropy radius, this substitution does not affect any of the properties in Section~\ref{sec:epo}.



\begin{figure}[!htb]
	\centering
	\begin{subfigure}[t]{.48\textwidth}
		\centering
		\includegraphics[width=1\textwidth]{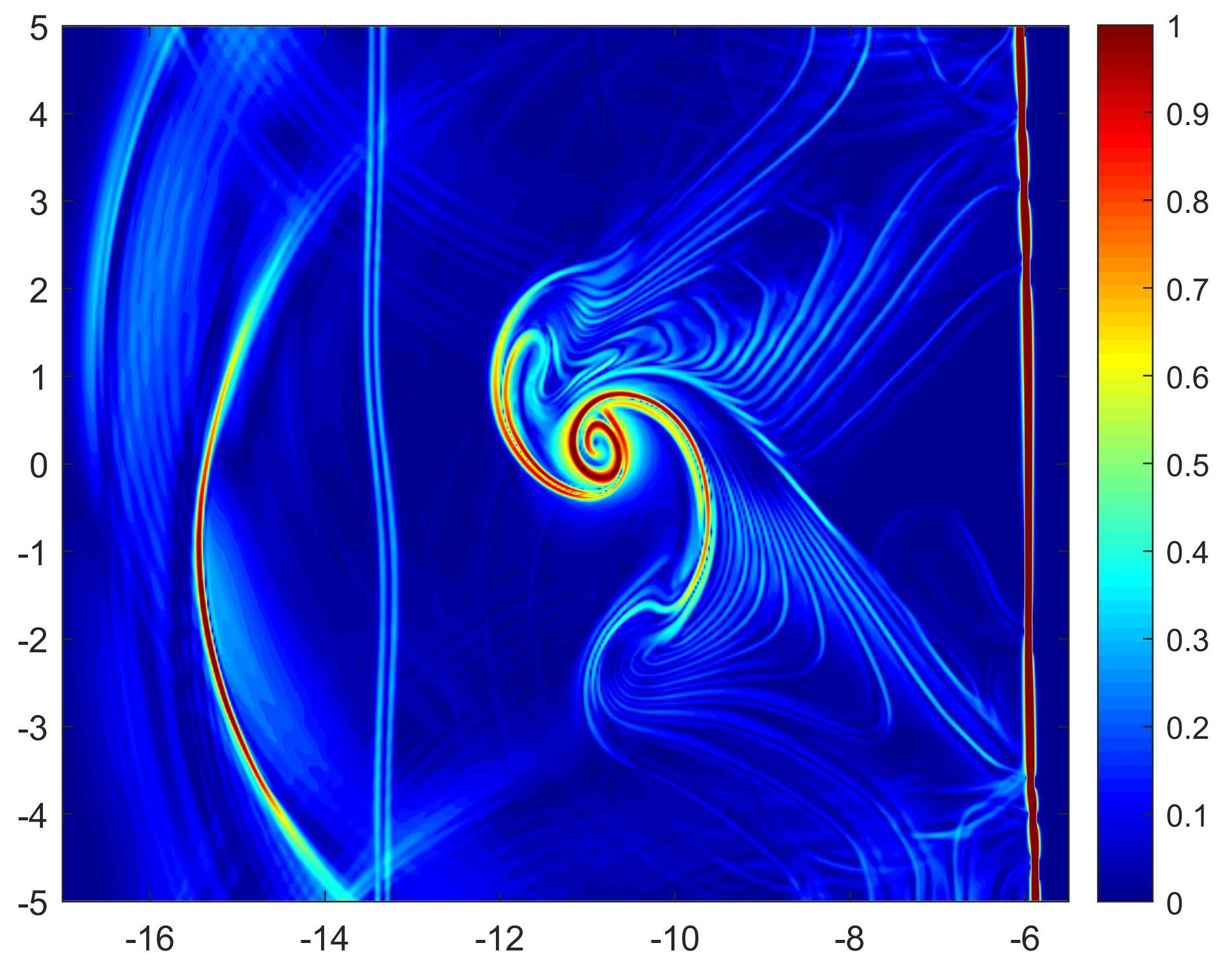}
		\caption{$\varepsilon_{\rm vtx} = 5$}
	\end{subfigure}
    \begin{subfigure}[t]{.48\textwidth}
		\centering
		\includegraphics[width=1\textwidth]{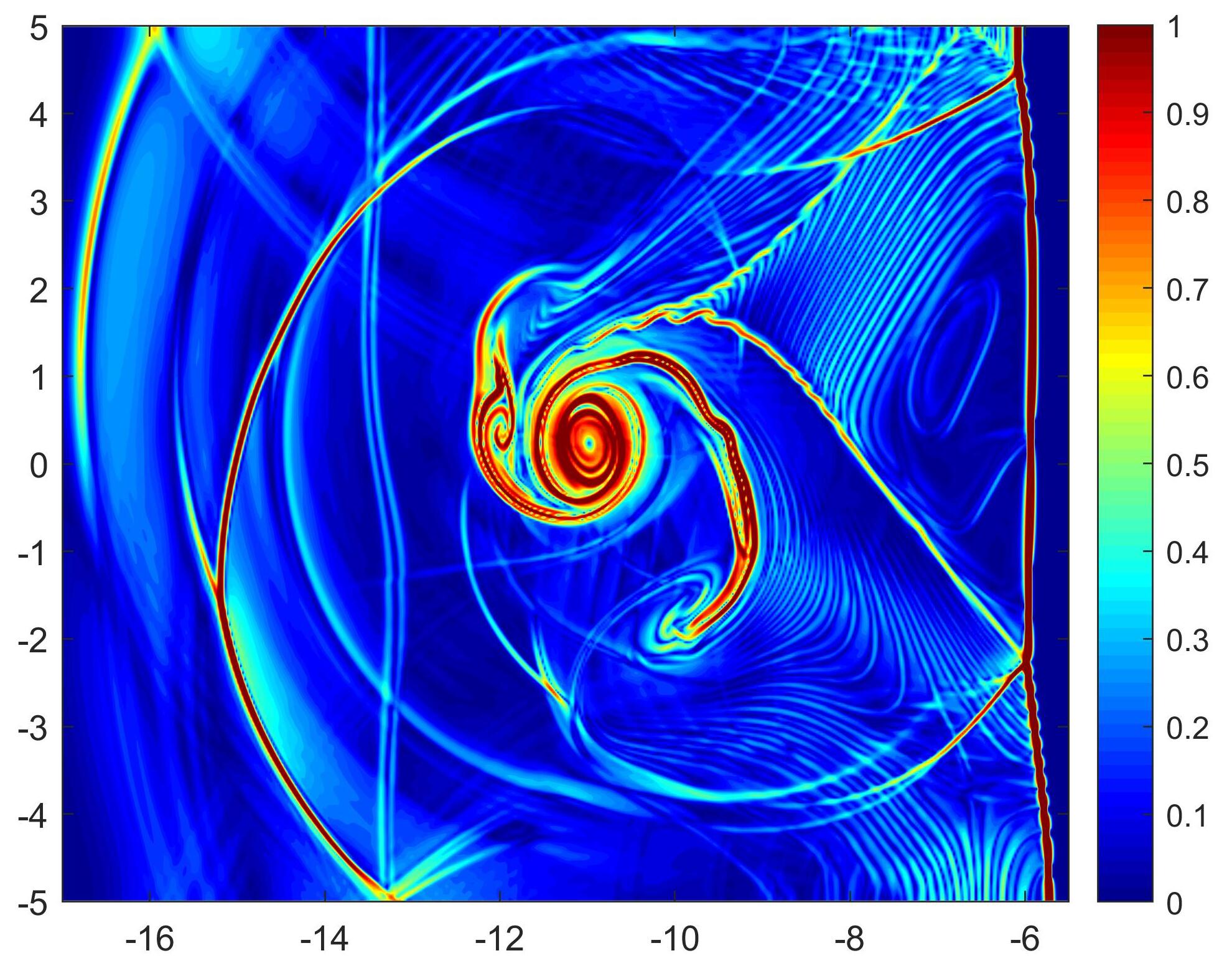}
		\caption{$\varepsilon_{\rm vtx} = 10.0828$}
	\end{subfigure}
	\caption{Test~\ref{Ex:SV}: Contour plots of ${\rm log}_{10}(1+|\nabla \rho|)$ at $t = 19$. Fifty equally distributed contour lines from 0 to 1 are displayed.}\label{Fig:2D_SV}
\end{figure}

\testcase{Shock--vortex interaction}{Ex:SV}
This test studies the interaction between a stationary shock and a moving vortex. A stationary shock is placed at $x=-6$, with the constant state $\mathbf{V}_{\rm L}$ behind it and the vortex superimposed on the region $x\ge-6$:
\begin{equation*}
    \mathbf{V}(x,y,0) =
    \begin{cases}
        \mathbf{V}_{\rm L}, & x < -6, \\[6pt]
        \mathbf{V}_{\rm vortex}(x,y), & x \ge -6,
    \end{cases}
\end{equation*}
with
\[
    \mathbf{V}_{\rm L} = (4.891497310766981, -0.388882958251919, 0, 11.894863258311670)^\top
\]
and the vortex state defined by
\[
    \mathbf{V}_{\rm vortex} =
    \begin{pmatrix}
        \bigl(\hat{\theta}(x,y)\bigr)^{\frac{1}{\Gamma-1}},\
        \dfrac{w_1(x,y)-v_c}{1-v_c\,w_1(x,y)},\
        \dfrac{\sqrt{1-v_c^2}\;w_2(x,y)}{1-v_c\,w_1(x,y)},\
        \bigl(\hat{\theta}(x,y)\bigr)^{\frac{\Gamma}{\Gamma-1}}
    \end{pmatrix}^\top.
\]
where
\begin{align*}
    \hat{\theta}(x,y) &= 1 - \frac{\varepsilon_{\rm vtx}^2}{28\pi^2}\,
        \exp\!\Bigl(1-\frac{x^2}{1-v_c^2}-y^2\Bigr), \\
    w_1(x,y) &= -y\,\hat{f}(x,y),\qquad
    w_2(x,y) = \frac{x}{\sqrt{1-v_c^2}}\,\hat{f}(x,y), \\
    \hat{f}(x,y) &= \sqrt{\frac{\beta}{1+\beta\bigl(\frac{x^2}{1-v_c^2}+y^2\bigr)}},\qquad
    \beta = \frac{\frac{\varepsilon_{\rm vtx}^2}{10\pi^2}\,
        \exp\!\bigl(1-\frac{x^2}{1-v_c^2}-y^2\bigr)}
        {1.8-\frac{\varepsilon_{\rm vtx}^2}{20\pi^2}\,
        \exp\!\bigl(1-\frac{x^2}{1-v_c^2}-y^2\bigr)}.
\end{align*}
The ideal EOS with $\Gamma = 1.4$ is used. The computational domain is $[-17,3]\times[-5,5]$ with reflective boundaries at $y=\pm5$, inflow at $x=3$ and outflow at $x=-17$. The simulation runs to $t=19$ on a uniform $800\times400$ mesh with the global COS module of Section~\ref{subsec:cos} active.

Two vortex strengths are considered: a mild vortex with $\varepsilon_{\rm vtx}=5.0$, and a severe one with $\varepsilon_{\rm vtx}=10.0828$, for which the pressure and the density at the vortex core are driven to near-vacuum values. In both cases, disabling the positivity module produces a negative density or pressure within a few time steps, leading to immediate code breakdown. With the full EPO projection, every nodal state remains admissible and the computation runs to the final time.

Figure~\ref{Fig:2D_SV} shows contour plots of $\log_{10}(1+|\nabla\rho|)$ at $t=19$. The deformed shock front, the primary vortex, and the train of secondary reflected waves behind it are all resolved, with sharp gradients and no visible
oscillations along the shock. The severe case differs from the mild one mainly in the strength of the shock deformation and the extent of the low-density core, both of which are captured without any special treatment.

\testcase{Relativistic jets}{Ex:jet}
Relativistic jets appear in many astrophysical settings, including active galactic nuclei and gamma-ray bursts. The numerical simulation of such flows is notoriously challenging due to the simultaneous presence of strong relativistic shocks, sharp shear layers, and Kelvin--Helmholtz instabilities, all of which can readily produce negative pressure or density and cause code failure without a rigorous bound-preserving mechanism.

We consider a pressure-matched, highly supersonic jet injected into a static ambient medium with state
\begin{align*}
    \mathbf{V}_{\text{a}}(x,y,0) &= \bigl(1,\,0,\,0,\,p_{\text{a}}\bigr)^\top,
\end{align*}
and three configurations of increasing severity:
\begin{enumerate}[label=(\arabic*),leftmargin=2em]
    \item $p_{\text{a}} = 2.3536240721718810\times10^{-5},\;
    v^b = 0.99,\; M_b = 50$;
    \item $p_{\text{a}} = 2.3966375079777595\times10^{-5},\;
    v^b = 0.999,\; M_b = 50$;
    \item $p_{\text{a}} = 2.3995344183272123\times10^{-7},\;
    v^b = 0.9999,\; M_b = 500$.
\end{enumerate}
The beam has rest-mass density $\rho^b = 0.1$ and enters through the inlet $\{y=0,\;|x| \le 0.5\}$ with vertical velocity $v^b$ and pressure matched to the ambient value. $M_b$ denotes the classical beam Mach number $v^b/c_s^b$. Exploiting the symmetry about $x=0$, we compute this example on $[0,12]\times[0,25]$ with $240\times500$ uniform cells. The displayed images are mirrored to $[-12,12]\times[0,25]$. The boundary conditions are reflective at $x=0$, beam inflow on the nozzle, and outflow on the remaining boundaries.

\begin{figure}[!htb]
	\centering
	\begin{subfigure}[t]{.32\textwidth}
		\centering
		\includegraphics[width=1\textwidth]{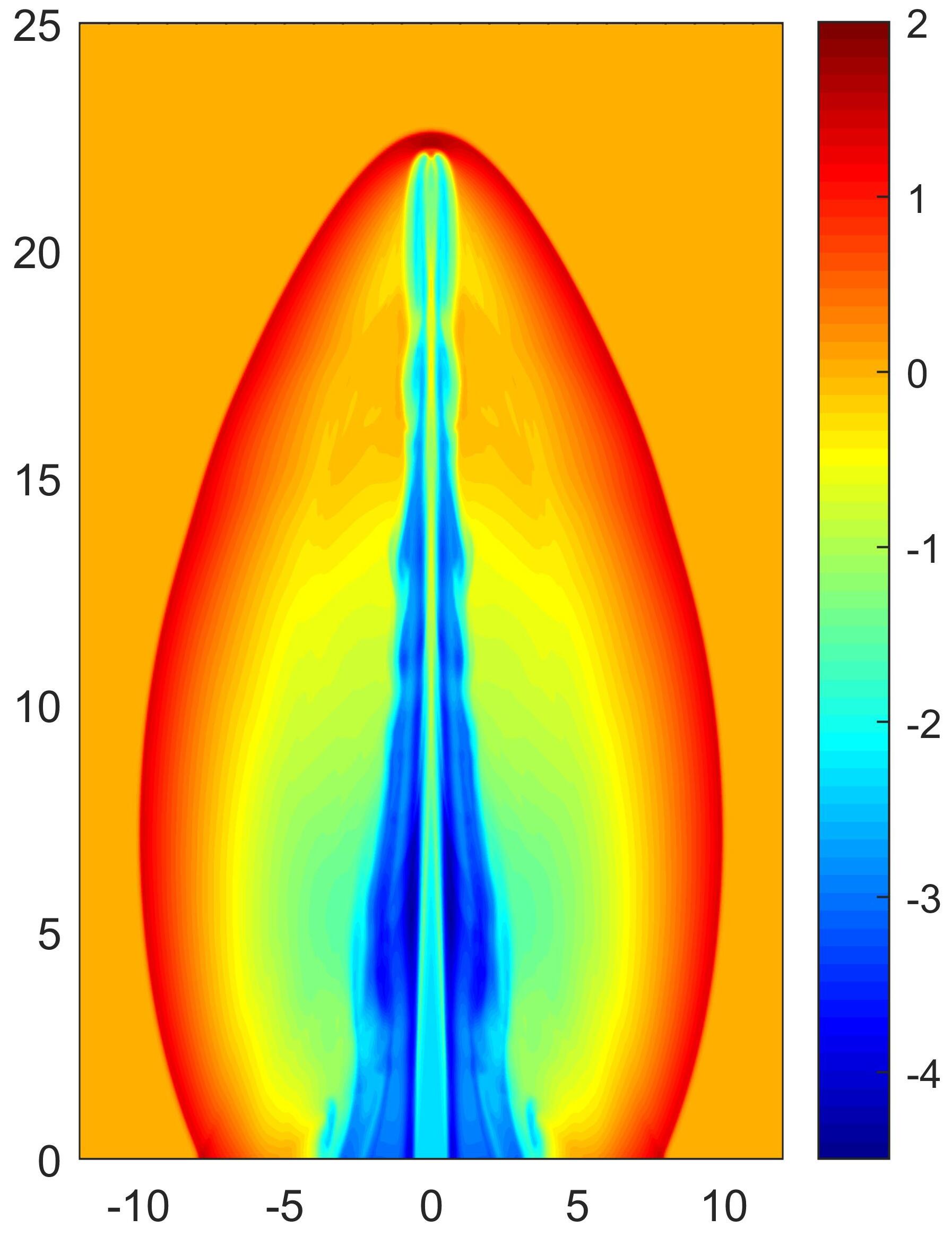}
	\end{subfigure}
    \begin{subfigure}[t]{.32\textwidth}
		\centering
		\includegraphics[width=1\textwidth]{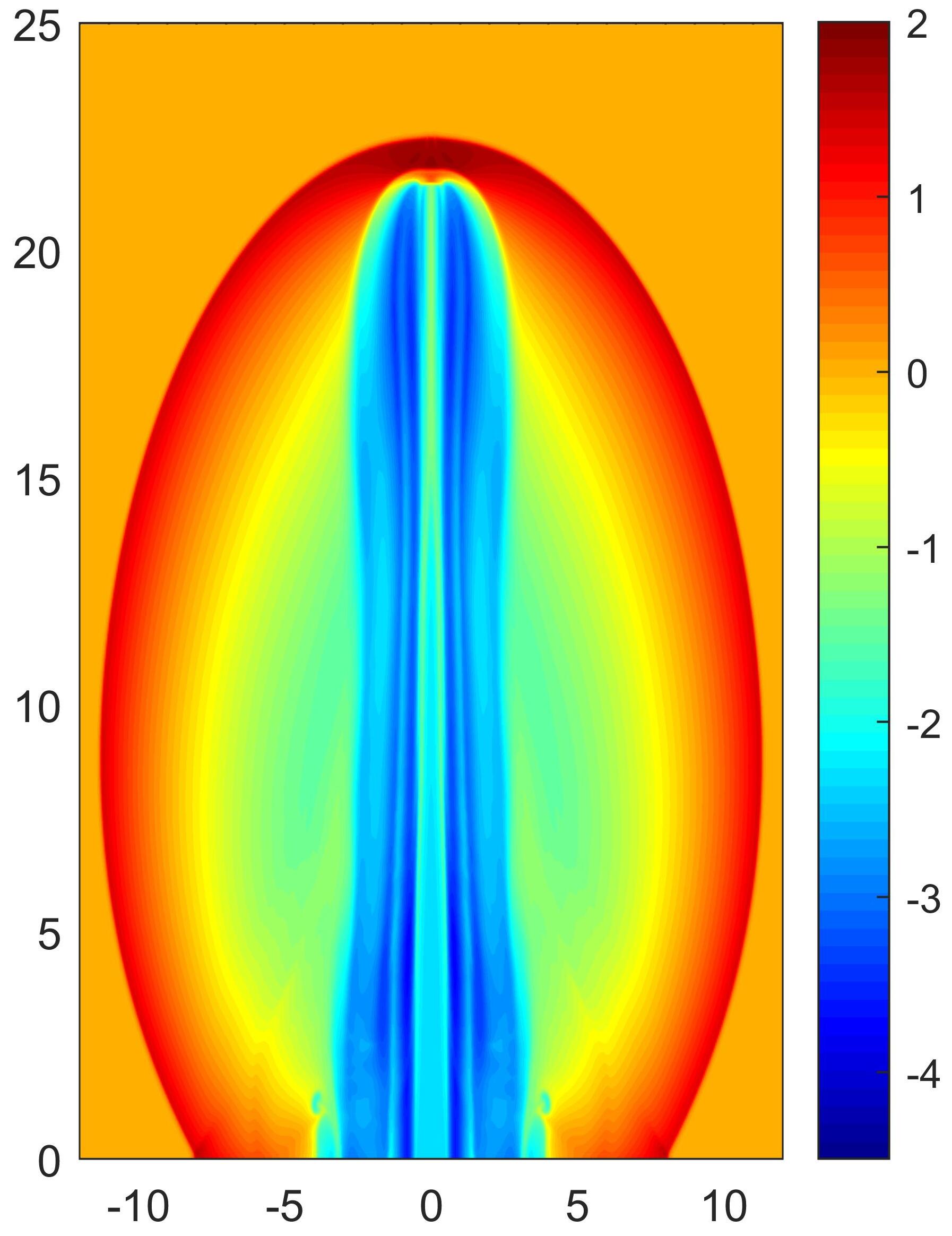}
	\end{subfigure}
 	\begin{subfigure}[t]{.32\textwidth}
		\centering
		\includegraphics[width=1\textwidth]{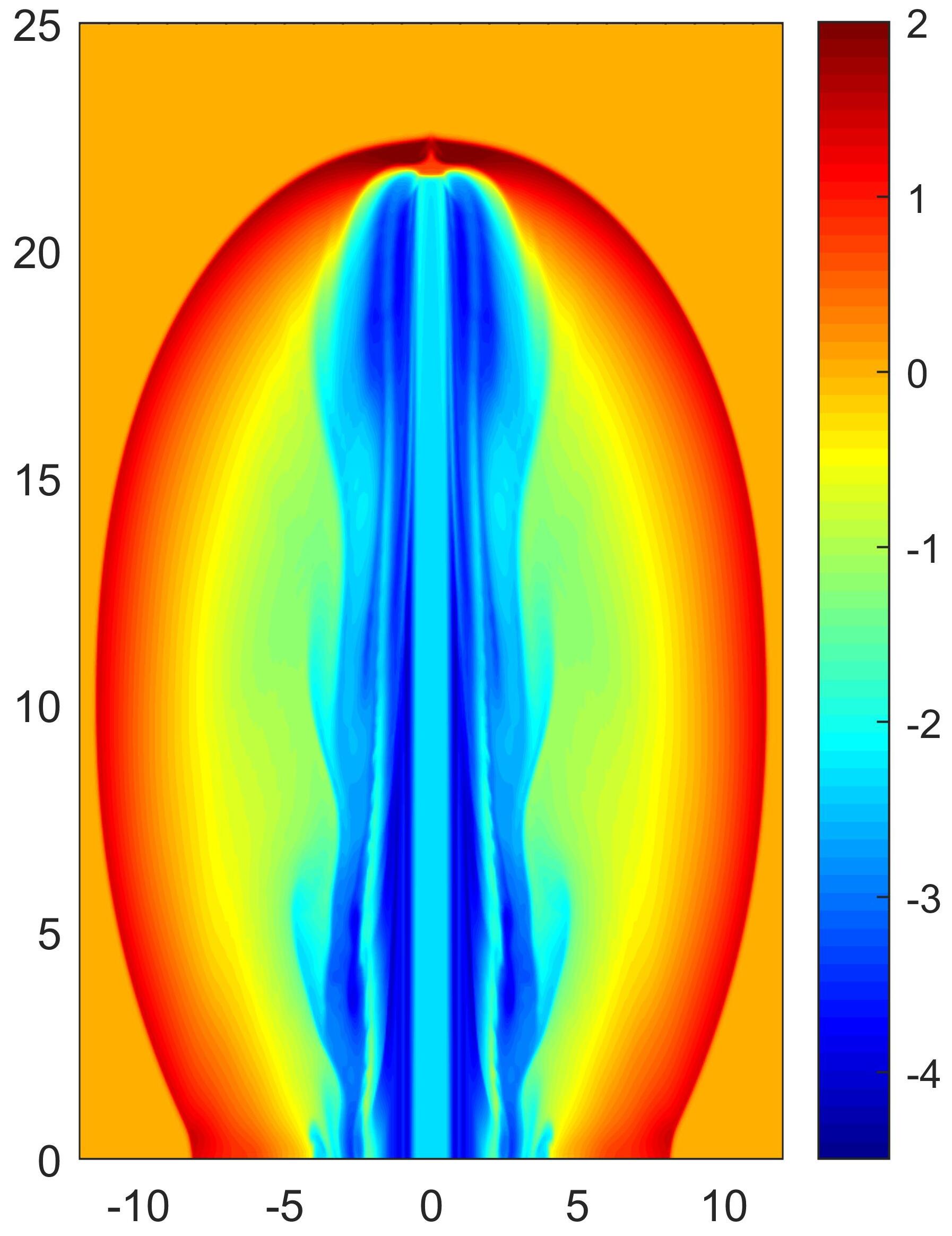}
	\end{subfigure}
	\caption{Test~\ref{Ex:jet}: Schlieren images of $\ln\rho$ obtained on $240\times500$ uniform cells with COS and displayed mirrored about $x=0$. From left to right: configurations (1) at $t=30$, (2) at $t=25$, and (3) at $t=23$.}\label{Fig:2D_Jet_lnrho}
\end{figure}

Figure~\ref{Fig:2D_Jet_lnrho} shows schlieren images of $\ln\rho$ for the three configurations. The leading bow shock, the cocoon, and the fine-scale Kelvin--Helmholtz roll-ups along the shear layer are resolved in all three regimes. The third configuration, whose beam Lorentz factor exceeds $70$, runs to $t=23$ without a single inadmissible state. Accuracy is not assessed in this test. We only note that the same $\mathbb{P}^2$ DG discretization fails during the initial transient in all three configurations if the positivity module is turned off.

\section{Conclusions}\label{sec:conclusion}

For the special relativistic Euler equations, physical admissibility and entropy stability cannot be arranged one after the other, because an inadmissible state has no entropy to correct, and the recovery through which every entropy is evaluated is well posed only on the admissible set. This paper has therefore treated them as one
coupled requirement and enforced both with a single cellwise projection. For any finite family of convex entropy pairs chosen by the user, the same conservative DG or finite volume candidate is corrected once, so that every quadrature state stays in the relativistic invariant domain and every member of the family satisfies its local strong entropy inequality, a global estimate followed by summation over the cells. Since SSP multistep integrators are convex combinations of forward Euler steps, the schemes are also fully discrete. Conservation, nodal admissibility, and the $M$ entropy inequalities consequently hold for the same computed state at every time step, whereas the entropy-stable RHD schemes we are aware of \cite{duan2019high,bhoriya2020entropy,duan2021entropy,wu2020entropy,duan2020high,xu2024high} give a semi-discrete estimate for a single pair and impose positivity, without which the entropy is not even defined, through a separate limiter.

What makes this possible is one structural fact. Relativistic causality bounds every directional wave speed by unity, so the Lax--Friedrichs average may be formed with $\alpha=1$ for every admissible state and every causal EOS. Being the average of the self-similar Riemann solution over the light cone, that average is itself admissible and inherits a two-point entropy inequality for each pair, since admissibility for the thermodynamic pair implies admissibility for all the others. A Gauss--Lobatto decomposition then transfers these inequalities to the high-order candidate average, and enforcing admissibility first guarantees that every primitive recovery and every entropy evaluation is performed on a state for which it is defined. The radial geometry and the SSP convexity argument are inherited from \cite{wu2026epo}, whereas the relativistic analysis that makes them applicable is developed here.

The computations bear this out. The design order is retained on smooth solutions whose troughs come within $10^{-5}$ of vacuum, and three entropy pairs in place of one remove the overshoot that persists behind the contact discontinuity, at the cost of two extra scalar entropy evaluations per trial radius. Moreover, no inadmissible nodal value occurs in any run, including interacting blast waves, shock heating at an incoming Lorentz factor of order $10^{5}$, a near-vacuum shock--vortex interaction, and jets with Lorentz factors above $70$, the last two of which break down within a few steps once the positivity module is switched off. Every reported entropy history decreases monotonically, in agreement with the theoretical estimates.

Limitations remain. The Gauss--Lobatto endpoint weight is only $\omega_1=1/56$ for $L=8$, which is a restrictive CFL condition, and sharper cell average decompositions such as \cite{cui2024optimal} should reduce this cost. The
multidimensional implementation, in addition, relies on line-separable Cartesian meshes, and an arbitrary finite entropy family is enforced at finite cost rather than the full infinite one. Beyond RHD, however, the construction uses the model only through the convexity of the admissible set, the light-cone bound, and the
generating condition for the entropy family, so the same route should be open for relativistic magnetohydrodynamics and for general relativistic hydrodynamics, where the primitive recovery is again implicit and the same ordering is forced.

% --------------------------------------------------------------------------------------------------

\section*{Data availability}
The data and code that support the findings of this study are available from the corresponding author upon reasonable request.

\section*{Declaration of competing interest}
The authors declare that they have no known competing financial interests or personal relationships that could have appeared to influence the work reported in this paper.

% --------------------------------------------------------------------------------------------------

\section*{Acknowledgement}
This work was partially supported by Science Challenge Project (No.~TZ2025007) and the Shenzhen Science and Technology Program (Grant Nos.~JCYJ20250604144300001 and RCJC20221008092757098).
% --------------------------------------------------------------------------------------------------

\appendix

\section{Primitive recovery and entropy evaluation}\label{app:primitive}

The entropy function is naturally expressed in primitive variables, but the numerical method evolves conservative variables. This appendix describes the primitive recovery procedure needed to evaluate the entropy profile in the limiter.

From the defining relations $D=\rho\gamma$, $\mathbf m=\rho h\gamma^2\mathbf v$ and $E=\rho h\gamma^2-p$, we obtain
\begin{equation*}
E+p=\rho h(\theta)\gamma^2,
\qquad
\mathbf m=(E+p)\,\mathbf v ,
\end{equation*}
so the recovery reduces to a single scalar unknown. Given
$\mathbf U=(D,\mathbf m,E)^\top\in\Gset$ and a trial pressure $p>0$, set
\begin{equation}\label{eq:recovery-chain}
\mathbf v(p):=\frac{\mathbf m}{E+p},
\qquad
\gamma(p):=\bigl(1-|\mathbf v(p)|^2\bigr)^{-1/2},
\qquad
\rho(p):=\frac{D}{\gamma(p)},
\qquad
\theta(p):=\frac{p}{\rho(p)},
\end{equation}
and define the residual
\begin{equation*}
\mathcal R(p)
:=
\rho(p)\,h\bigl(\theta(p)\bigr)\,\gamma(p)^2-\bigl(E+p\bigr).
\end{equation*}
The pressure of the state $\mathbf U$ is the root of $\mathcal R$, and the remaining primitive variables then follow from \eqref{eq:recovery-chain}.
For admissible conservative data, the existence and uniqueness of the primitive variables are standard, and we refer to \cite{WuTang2015,WuTang2017ApJS,cai2024provably}. Recently, a provably convergent, constraint-preserving Newton method for primitive recovery was introduced in \cite{cai2024provably, yue2026robust}. The positivity module ensures that every nodal state entering an entropy profile lies in $\Geps$, so the recovery map stays uniformly bounded away from vacuum and from $|\mathbf v|=1$, and the safeguarded iteration converges in a fixed small number of steps.

For the canonical entropy pair, once $(\rho,\theta,\mathbf v)$ are recovered,
\[
\E(\mathbf U)=-D\left(-\ln\rho+\int_1^\theta \frac{e'(\xi)}{\xi}\,\dd\xi\right).
\]
For the specific EOS examples in \ref{app:eos}, the entropy integral is explicit, and no quadrature is needed.

\section{Canonical entropy pairs for standard equations of state}\label{app:eos}

This appendix gives the canonical specific entropy $S(\mathbf U)$ and the associated entropy pair
\[
\E(\mathbf U)=-DS(\mathbf U),
\qquad
\Qflux_i(\mathbf U)=-Dv_iS(\mathbf U).
\]

\subsection{Ideal-gas EOS}

For the ideal EOS
\[
h(\theta)=1+\frac{\Gamma}{\Gamma-1}\theta,
\qquad 1<\Gamma\le 2,
\]
one has
\[
e(\theta)=\frac{\theta}{\Gamma-1},
\qquad
S(\mathbf U)=\frac{1}{\Gamma-1}\ln\frac{p}{\rho^\Gamma}.
\]
Hence
\[
\E_{ID}(\mathbf U)=-\frac{D}{\Gamma-1}\ln\frac{p}{\rho^\Gamma},
\qquad
\Qflux_{ID,i}(\mathbf U)=-\frac{Dv_i}{\Gamma-1}\ln\frac{p}{\rho^\Gamma}.
\]

\subsection{RC EOS}

For the RC EOS
\[
h(\theta)=\frac{2(6\theta^2+4\theta+1)}{3\theta+2},
\]
the specific entropy is
\[
S_{RC}(\mathbf U)
=
-\ln\rho
+\frac32\ln\theta
+\frac32\ln(3\theta+2)
-\frac{3}{3\theta+2}
+C_{RC},
\]
with $C_{RC}=\frac35-\frac32\ln5$. Therefore
\[
\E_{RC}(\mathbf U)=-D S_{RC}(\mathbf U),
\qquad
\Qflux_{RC,i}(\mathbf U)=-Dv_i S_{RC}(\mathbf U).
\]

\subsection{IP EOS}

For the IP EOS
\[
h(\theta)=2\theta+\sqrt{1+4\theta^2},
\]
one finds
\[
S_{IP}(\mathbf U)
=
-\ln\rho
+\ln\theta
+2\ln\bigl(2\theta+\sqrt{1+4\theta^2}\bigr)
+C_{IP},
\]
with $C_{IP}=-2\ln(2+\sqrt5)$. Therefore
\[
\E_{IP}(\mathbf U)=-D S_{IP}(\mathbf U),
\qquad
\Qflux_{IP,i}(\mathbf U)=-Dv_i S_{IP}(\mathbf U).
\]

\subsection{TM EOS}

For the TM EOS
\[
h(\theta)=\frac52\theta+\sqrt{1+\frac94\theta^2},
\]
the specific entropy is
\[
S_{TM}(\mathbf U)
=
-\ln\rho
+\frac32\ln\theta
+\frac32\ln\left(\frac32\theta+\sqrt{1+\frac94\theta^2}\right)
+C_{TM},
\]
with $C_{TM}=-\frac32\ln\left(\frac32+\sqrt{\frac{13}{4}}\right)$. Therefore
\[
\E_{TM}(\mathbf U)=-D S_{TM}(\mathbf U),
\qquad
\Qflux_{TM,i}(\mathbf U)=-Dv_i S_{TM}(\mathbf U).
\]

\bibliographystyle{elsarticle-num}
\bibliography{rhd_epo_refs}

\end{document}